\documentclass[11pt]{amsart}

\usepackage{amsmath,amsthm, amscd, amssymb, amsfonts,mathtools}
\usepackage[all]{xy}
\usepackage{enumitem}

\usepackage{t1enc}
\usepackage[latin1]{inputenc}

\usepackage{fontsmpl}
\usepackage{srcltx}

\allowdisplaybreaks

\usepackage{color}

\numberwithin{equation}{section}\theoremstyle{plain}

\newtheorem{theorem}{Theorem}[section]
\newtheorem*{maintheorem}{Main Theorem}

\newtheorem{lem}[theorem]{Lemma}
\newtheorem{cor}[theorem]{Corollary}
\newtheorem{pro}[theorem]{Proposition}

\theoremstyle{definition}
\newtheorem{definition}[theorem]{Definition}
\newtheorem{exa}[theorem]{Example}

\newtheorem{conjecture}{Conjecture}

\theoremstyle{remark}
\newtheorem{rem}[theorem]{Remark}

\newcommand{\ydh}{{}^{H}_{H}\mathcal{YD}}

\newcommand{\ydf}{{}^{F}_{F}\mathcal{YD}}

\def\G{\mathbb{G}}

\newcommand\id{\operatorname{id}}

\newcommand\Alg{\operatorname{Alg}}
\newcommand\Cleft{\operatorname{Cleft}}

\newcommand\gr{\operatorname{gr}}

\newcommand\Sh{\operatorname{Sh}}

\newcommand{\uvi}{\mathtt u}

\def\q{\mathbf{q}}

\def\k{\Bbbk}

\def\ot{\otimes}

\def\Z{\mathbb{Z}}
\def\N{\mathbb{N}}

\def\B{\mathfrak{B}}
\def\L{\mathfrak{L}}
\def\z{\mathfrak{z}}
\def\n{\mathfrak{n}}

\def\eps{\epsilon}
\def\mT{\mathcal{T}}

\def\mL{\mathcal{L}}
\def\mP{\mathcal{P}}
\def\mO{\mathcal{O}}
\def\mR{\mathcal{R}}

\def\mA{\mathcal{A}}
\def\mE{\mathcal{E}}
\def\mH{\mathcal{H}}

\def\u{\mathfrak{u}}

\def\bx{\mathbf{x}}
\def\bg{\mathbf{g}}
\def\bchi{\bs\chi}
\def\bbeta{\bs\beta}

\def\ta{\texttt{a}}
\def\tb{\texttt{b}}

\newcommand{\J}{{\mathcal J}}
\newcommand{\Gc}{{\mathcal G}}
\newcommand{\D}{{\mathcal D}}
\newcommand\I{\mathbb I}

\newcommand{\cf}{{\it cf.~}}

\newcommand{\loc}{{\it loc.cit.~}}
\def\bs{\boldsymbol}

\def\qb{\mathfrak{q}}

\def\lg{\langle}
\def\rg{\rangle}

\def\pf{\begin{proof}}
\def\epf{\end{proof}}

\renewcommand{\a}[1]{^{#1}}

\begin{document}


 \title[Pointed Hopf algebras of type $G_2$]{\small  Liftings of Nichols algebras of diagonal type \\ III. Cartan type $G_2$}
\author[Garc\'ia Iglesias; Jury Giraldi]{Agust\'in Garc\'ia Iglesias; Jo\~ao Matheus Jury Giraldi}

\address{\noindent A.G.I.: FaMAF-CIEM (CONICET), Universidad Nacional de C\'ordoba,
Me\-dina A\-llen\-de s/n, Ciudad Universitaria, 5000 C\' ordoba, Rep\'
ublica Argentina.}

\address{\noindent J.M.J.G.: Instituto de Matem\'atica, 
Universidade Federal do Rio Grande do Sul,
Rio Grande do Sul, Brazil.}

\email{aigarcia@famaf.unc.edu.ar, joaomjg@gmail.com}

\thanks{The work was partially supported by CONICET,
FONCyT-ANPCyT, Secyt (UNC), the MathAmSud project
GR2HOPF, CNPq (Brazil)}

\begin{abstract}
We complete	 the classification of Hopf algebras whose infinitesimal braiding is a principal Yetter-Drinfeld realization of a braided vector space of Cartan type $G_2$ over a cosemisimple Hopf algebra.

We develop a general formula for a class of liftings in which the quantum Serre relations hold. We give a detailed explanation of the procedure for finding the relations, based on the recent work of Andruskiewitsch, Angiono and Rossi Bertone.
\end{abstract}
\maketitle
\section{Introduction}

This article belongs to the series initiated in \cite{AAG}, and followed by \cite{AGI2}, with the purpose of computing all liftings of braided vector spaces of diagonal type $(V,c)$ whose Nichols algebra $\B(V)$ is finite-dimensi\-onal. 

In this part we focus on:
\begin{itemize} [leftmargin=*]\renewcommand{\labelitemi}{$\circ$}
\item $(V,c)$ of Cartan type $G_2$. That is, we assume there are $N>3$, $q$ a primitive $N$-th root of 1 and a basis $\{x_1,x_2\}$ of $V$ such that the braiding is determined by a matrix $\qb=\left(\begin{smallmatrix}
q&q_{12}\\ q_{21}&q^3
\end{smallmatrix}\right)\in\k^{2\times 2}$ with $q_{12}q_{21}=q^{-3}$ as:
\begin{align}\label{eqn:braiding-G2}
\begin{split}
c(x_1\ot x_1)&=q x_1\ot x_1; \qquad c(x_1\ot x_2)=q_{12} x_2\ot x_1; \\
c(x_2\ot x_1)&=q_{21} x_1\ot x_2; \quad c(x_2\ot x_2)=q^3 x_2\ot x_2.
\end{split}
\end{align}
\item $H$ a cosemisimple Hopf algebra with a principal realization $(\chi_i,g_i)_{i=1,2}$, so that $V\in \ydh$. That is, $g_1,g_2\in Z(G(H))$ and $\chi_1,\chi_2\in\Alg(H,\k)$ are such that $\chi_i(g_j)=q_{ji}$, $i,j\in\{1,2\}$.
\end{itemize}

We recall that a Hopf algebra $L$ is called a lifting of $V\in\ydh$ when $\gr L\simeq \B(V)\# H$. We refer the reader to the article \cite{AAG}, the first of the series, based on \cite{AAGMV}, for a description of the program for computing all liftings of braided vector spaces of diagonal type with a realization $V\in\ydh$. 

Set $B=\B(V)\# H$. In a sentence, the program consists in constructing a subset $\Cleft'(B)\subseteq \Cleft(B)$ of right cleft objects for $B$ in such a way that for every $X\in \Cleft'(B)$ the left Schauenburg Hopf algebra $A=L(X,B)$, see \cite{S}, satisfies $\gr A\simeq B$ and checking that indeed, every lifting can be obtained in this way. 

Assume that $(N,3)=1$, $N>7$ and that $H=\k\Gamma$, $\Gamma$ a finite abelian group. A family of Hopf algebras $\u(\D,\mu)$ was defined in \cite{AS-annals}, as a particular case of a more general construction, such that $\gr \u(\D,\mu)\simeq \B(V)\# H$. Furthermore, \cite[Classification Theorem 0.1]{AS-annals} shows that
when $(N,210)=1$ and $A$ is a lifting, then there is $(\D,\mu)$ such that $A\simeq \u(\D,\mu)$.

In this article we give a complete presentation of every lifting of Cartan type $G_2$, without restrictions on the parameters and for any cosemisimple Hopf algebra $H$. We obtain new examples of Hopf algebras $\u(\bs\mu)$ and $\u(\bs\lambda,\bs\mu)$, see Definitions \ref{def:generic} and \ref{def:deformed}. 

The case in which the quantum Serre relations can be deformed has a very particular nature. This phenomenon only occurs when the braiding matrix is $\qb=\qb_d$ (hence $N=7$) as in \eqref{eqn:quantum-def}.
As in \cite{AAG}, see also \cite{helbig}, in this case the relations are deformed in positive terms of the coradical filtration. It presents, however, a new issue: the description of the algebras $\u(\bs\lambda,\bs\mu)$ in Definition \ref{def:deformed} takes us {\it more than 20 pages}, and it is attained only after a tedious combinatorial puzzle, whose solution involves many running days of a computer program. We see no symmetry in the presentation of these algebras, rather we observe a chaotic outburst of coefficients preceding a big number of terms of PBW generators. A similar, although shorter, situation is already present in the deformations in \cite{AAG} and \cite{helbig}. This case looks like a singularity that explodes when one deforms the quantum Serre relations, responsible for the nice behavior of the power of the root vectors, a role that becomes evident once we attempt to tamper with them. 

The computation of the liftings thus reveals itself as a full challenge, in the sense that we must not only deal with the question on how to solve this problem, but also on how to communicate the results. 

On the other hand, we point out that we obtain in Definition \ref{def:generic} a closed and concise description of the generic case, which extends the results in \cite{AS-annals} mentioned above. We give a full description of the relations, with no recurrence formula involved. For this, we use a technique introduced in \cite{AAR2} to compute the coproducts of the powers  of root vectors.

\smallbreak

Our main contribution reads as follows. 
\begin{maintheorem}
Let $H$ be a cosemisimple Hopf algebra and let $V$ be a braided vector space of diagonal type Cartan $G_2$ with a principal Yetter-Drinfeld realization $V\in\ydh$; let $\qb$ be the braiding matrix. 

Let $L$ be a Hopf algebra such that $L_0\simeq H$ and its infinitesimal braiding is given by $V\in\ydh$. Then
\begin{enumerate}[leftmargin=*]
\item If $\qb\neq \qb_d$, then there is $\bs\mu\in\k^{\Delta_+}$ such that $L\simeq \u(\bs\mu)$.
\item If $\qb=\qb_d$, then there are $\bs\lambda\in\k^2$ and $\bs\mu\in\k^{\Delta_+}$ such that $L\simeq \u(\bs\lambda,\bs\mu)$.
\end{enumerate}
\end{maintheorem}
\pf
Recall that every such $L$ is a lifting of $V$ by \cite[Theorem 2]{A-nichols}.
Now, (1) is Theorem \ref{thm:2steps}; (2) is Theorem \ref{thm:deformed}.
\epf
The isomorphism classes of the algebras $\u(\bs\mu)$, resp.~$\u(\bs\lambda,\bs\mu)$, are characterized in terms of certain symmetries of the families  $\bs\mu$, resp.~$(\bs\lambda,\bs\mu)$; see \cite[Theorem 3.9]{AAG}. The more general class of braided vector spaces of {\it standard} type $G_2$ is currently under investigation.

As a byproduct of our research, we present in \S \ref{sec:two-steps} a recursive method to compute all liftings in which the quantum Serre relations still hold.

The paper is organized as follows. In Section \ref{sec:prelim} we fix the notation and collect some standard facts about Nichols algebras and their deformations. We also recall the basics of our strategy to find all of these liftings. In Section 
\ref{sec:two-steps} we apply this procedure to give a generic presentation of every lifting satisfying some technical conditions. In Section \ref{sec:cartan} we prove our main results, that is the explicit description of every lifting of a braided vector space of diagonal type $G_2$. We  include an Appendix with some comments on the use of the computer to perform some of the computations.

\section*{Acknowledgments}
We thank Iv\'an Angiono for sharing his knowledge about Nichols algebras with us; this article could not have been completed without his generous assistance. We  thank Fiorela Rossi Bertone for interesting conversations and Nicol\'as Andruskiewitsch for his constant support and thoughtful advice. We also thank Fernando Fantino for his assistance with \texttt{GAP}. We thank the referee for his/her remarks, which  improved the presentation of this article.

\section{Preliminaries}\label{sec:prelim}
We work over an algebraically closed field of characteristic zero $\k$. If $A$ is a $\k$-algebra and $S\subseteq A$ is a subset, then we denote by $\lg S\rg$ the ideal generated by $S$ in $A$. We write $A^\times$ for the group of invertible elements in $A$ and $\Alg(A,\k)$ for the set of algebra maps $A\to \k$. If $G$ is a group, then $Z(G)$ denotes its center. If $S\subseteq G$ is a subset, then $\lg S\rg$ denotes the subgroup generated by $S$. If $n\in\N$, we write $\I_n=\{1,\dots,n\}$. For $n\geq 2$, we set $\G_n$ for the group of $n$-th roots of unity and fix $\G_n'\subset \G_n$ the subset of primitive roots. 

If $H$ is a Hopf algebra, then we denote by $G(H)$ the group of group-like elements of $H$ and by $\mP(H)$ the subspace of its primitive elements. We denote by $\ydh$ the category of Yetter-Drinfeld modules over $H$. If $M\in\ydh$ and $\delta$, resp. $\cdot$, denotes the $H$-coaction, resp. the $H$-action, on $M$, then we set, for each $g\in G(H)$, $\chi\in\Alg(H,\k)$:
\begin{align*}
M_g&=\{x\in M| \delta(x)=g\ot x\}, & M^\chi=\{x\in M| h\cdot x=\chi(h)x, \forall\, h\in H\}.
\end{align*} 
In addition, we set $M_g\a{\chi}=M_g\cap M\a{\chi}$.

\subsection{Nichols algebras of diagonal type}

We review the basics on braided vector spaces $(V,c)$ of diagonal type, their Nichols and distinguished pre-Nichols algebras $\B(V)$ and $\widetilde{\B}(V)$, and the Lusztig algebras $\L(V)$.

\subsubsection{Braided vector spaces of diagonal type}
Let $(V,c)$ be a  braided vector space of \emph{diagonal type}, with $\dim V=\theta$. That is, there is a basis $(x_i)_{i\in\I_\theta}$ and a {\it braiding matrix} $\q=(q_{ij})_{i,j\in\I_\theta}$ in such a way that 
\[
c(x_i\ot x_j)=q_{ij}x_j\ot x_i.
\] 
We fix some notation, we refer the reader to \cite{AAG} and references therein for unexplained terminology. Let  $\left(g_i, \chi_i\right)_{i \in \I_\theta}$ be a principal realization of $V$ over $H$. This is a collection of group-like elements $g_i\in Z(G(H))$ and maps $\chi_j\in\Alg(H,\k)$ satisfying $\chi_j(g_i)=q_{ij}$, $i,j\in\I_\theta$. In this way $V$ becomes an object in $\ydh$, so that the braiding coincides with the categorical one, by declaring $x_i\in V^{\chi_i}_{g_i}$, $i\in\I_\theta$. In particular, we set
\begin{align}\label{eqn:Gamma}
\Gamma &= \lg g_1,\dots,g_\theta\rg < Z(G(H)).
\end{align}
Let $\J(V)\subset T(V)$ be the defining ideal of the Nichols algebras $\B(V)$ and fix $\Gc$ the set of generators of $\J(V)$ as computed in \cite{A-nichols}, for every case in the list of \cite{H-classif}.

Let $\Delta=\Delta^{\q}$ denote the {\it generalized root system} associated to $\q$. We write $\Delta_+=\{\beta_1,\dots,\beta_M\}$ for the subset of positive roots; this is equipped with a convex order $\beta_M>\dots>\beta_1$. Let $\alpha_i$, $i\in\I$, be the simple roots. We denote by $x_\alpha$, $\alpha\in\Delta_+$ the corresponding {\it root vectors}; hence $x_{\alpha_i}=x_i$, $i\in\I_\theta$. If $\alpha=n_1\alpha_{i_1}+\dots+n_k\alpha_{i_k}\in\Delta_+$, $k,n_1\dots, n_k\in\N$, then we set
\begin{align*}
g_{\alpha}&:=g_{i_1}^{n_1}\dots g_{i_k}^{n_k}\in\Gamma, &  \chi_{\alpha}&:=\chi_{i_1}^{n_1}\dots  \chi_{i_k}^{n_k}\in \Alg(H,\k).
\end{align*}

\subsubsection{Pre-Nichols and Lusztig algebras}\label{sec:normal}

We consider the subset $\mO=\mO(V)\subseteq\Delta_+$ of {\it Cartan roots} \cf \cite[(20)]{A-distinguished}. In particular, the {\it distinguished pre-Nichols algebra} $\widetilde{\B}(V)$ from \cite{A-distinguished} is defined as the  quotient of $T(V)$ modulo the ideal generated by $\Gc\setminus \{x_\alpha^{N_\alpha}|\alpha\in \mO(V)\}$ -- and possibly adding some redundant relations for $\B(V)$. 
We set $Z(V)=\k\lg x_\alpha^{N_\alpha}|\alpha\in \mO(V)\rg \subset \widetilde{\B}(V)$. A key fact is \cite[Theorem 4.13]{A-distinguished}, which in this context reads as follows:
\begin{theorem} \label{thm:normal}
$Z(V)\subset \widetilde{\B}(V)$ is a normal Hopf subalgebra.\qed
\end{theorem}
In particular, this determines an exact sequence
\begin{align*}
Z(V)\stackrel{\iota}{\longrightarrow} \widetilde{\B}(V)\stackrel{\pi}{\longrightarrow}\B(V).
\end{align*}
We write $\z(V)$ for the graded dual of $Z(V)$. In turn, we denote by $\L(V)$ the {\it Lusztig algebra} of $V$: this is defined in \cite{AAR1} as the graded dual of $\widetilde{\B}(V)$ and it fits into an extension of braided Hopf algebras 
\begin{align}\label{eqn:sequence}
\B(V)\stackrel{\pi^*}{\longrightarrow}\L(V)\stackrel{\iota^*}{\longrightarrow}\z(V),
\end{align}
see \cite[Proposition 3.2]{AAR2}. We remark that the convolution product on these dual Hopf algebras is defined as
\begin{align}\label{eqn:convolution}
\lg f*g, x\rg=\lg g, x_{(1)}\rg\lg f, x_{(2)}\rg.
\end{align}

We recall that $\L(V)$ is generated as an algebra by the {\it divided powers} $y_\beta^{(n)}$, $\beta\in\Delta_+$, $n\in\N$. Here $y_\beta^{(n)}\in \L(V)\simeq \widetilde{\B}(V)^*$ is such that 
\begin{align}\label{eqn:divided-powers}
\lg y_\beta^{(n)},x_\alpha^m\rg=\delta_{\alpha,\beta}\delta_{n,m}, \quad \alpha,\beta\in\Delta_+,\,n,m\in\N. 
\end{align}

\subsection{Twist equivalent braided vector spaces}\label{sec:twist-equiv}

Let $(V,c)$ and $(\hat{V},\hat{c})$ be two  braided vector spaces of diagonal type of dimension $\theta\in\N$. 

Let $\q=(q_{ij})_{i,j\in\I_\theta}$ and $\hat{\q}=(\hat{q}_{ij})_{i,j\in\I_\theta}$ denote the corresponding braiding matrices. We recall that  $V$ and $\hat{V}$ are {\it twist equivalent} if and only if
\begin{align*}
q_{ij}q_{ji}=\hat{q}_{ij}\hat{q}_{ji} \ \text{ and } \ q_{ii}=\hat{q}_{ii}, \qquad i,j\in \I_\theta.
\end{align*}
In particular, $\Delta^{\q}=\Delta^{\hat{\q}}=\{\beta_1,\dots,\beta_M\}$. 
We write $(x_i)_{i\in\I_\theta}$ and $(\hat{x}_i)_{i\in\I_\theta}$ for the bases of $V$ and $\hat{V}$, respectively.

In this setting, the Nichols algebras $\B(\hat{V})$ and $\B(V)$ are closely related, see Proposition \ref{pro:nichols-twist}. We need some notation to explain this relation.

We fix $F$ the free group on $\theta$ generators $g_1,\dots,g_\theta$; so there are principal realizations $(g_i,\chi_i)_{i\in\I_\theta}$, resp. $(g_i,\hat{\chi}_i)_{i\in\I_\theta}$, of $V\in\ydf$, resp. $\hat{V}\in\ydf$. We consider the group cocycle $\sigma:F\times F\to \k^\times$ defined by 
\[
\sigma(g_i,g_j)=\begin{cases}
\hat{q}_{ij}q_{ij}^{-1}, & i\leq j\in\I_\theta;\\
1, & i\geq j\in\I_\theta.
\end{cases}
\]
We fix a family of scalars $(t_\alpha)_{\alpha\in\Delta}$ as defined recursively in \cite[\S 4.1]{A-convex} via
\begin{align}\label{eqn:t-alpha}
t_{\alpha_i}&=1, \quad i\in\I_\theta; & t_\alpha &=\sigma(\alpha',\alpha'')t_{\alpha'} t_{\alpha''},
\end{align} 
where $\alpha\in\Delta_+$ is not simple and $(l_{\alpha'},l_{\alpha''})\coloneqq\Sh(l_\alpha)$ is the {\it Shirshov decomposition} of the Lyndon word $l_\alpha$ associated to the root $\alpha$; see \cite[Definition 2.2, Corollary 3.7]{A-convex} for details.

For every $(n_i)_{i\in\I_M}\in \Z_{\geq 0}\a{M}$, we set, following \cite[(4.6)]{A-convex}:
\[
f(\hat{x}_{\beta_M}^{n_M}\dots \hat{x}_{\beta_1}^{n_1})=\prod_{i<j\in\I_\theta}\sigma(\beta_j,\beta_i)^{n_in_j}\prod_{i\in\I_\theta}\sigma(\beta_i,\beta_i)^{\binom{n_i}{2}}t_{\beta_i}^n.
\]
If $\hat{x}\in \B(\hat{V})_h$, $h\in F$, then there are $h_i,h^i\in F$, $i\in I$, with $h_ih^i=h$, with
 \[
 \Delta(\hat{x})=\sum_{i\in I} \hat{x}_i\ot \hat{x}^i, \qquad \text{and }\hat{x}_i\in \B(V)_{h_i},  \hat{x}^i\in \B(V)_{h^i}.
 \]
  We consider the map $\Delta_\sigma:\B(\hat{V})\to \B(\hat{V})\ot\B(\hat{V})$ defined via
\begin{align} \label{eqn:coproddeformed}
\Delta_\sigma(\hat{x})= \sum_{i\in I}\sigma(h_i,h^i)\hat{x}_i\ot \hat{x}^i, \qquad \hat{x}\in \B(\hat{V}).
\end{align}
We recall that this map is defined in \cite{AS-pointed} by considering a right $\k\Gamma\ot \k\Gamma$-action $\leftharpoonup$ on $\B(\hat{V})\ot\B(\hat{V})$, so $\Delta_\sigma(\hat{x})=\Delta(\hat{x})\leftharpoonup \sigma$.
We use a Sweedler's type notation to write $\Delta_\sigma(\hat{x})=\sigma(h_{(1)},h_{(2)})\hat{x}_{(1)}\ot \hat{x}_{(2)}$. 

The next proposition establishes a twisted equivalence between the braided Hopf algebras $\B(V)$ and $\B(\hat{V})$. 
As a byproduct, it presents a formula to find an expression of the coproducts of the powers of root vectors $(x_\alpha^{N_\alpha})_{\alpha\in\Delta_+}$ in terms of the coproduct of the powers of root vectors $(\hat{x}_\alpha^{N_\alpha})_{\alpha\in\Delta_+}$ coming from a twist equivalent braided vector space.

\begin{pro}\cite[Proposition 3.9]{AS-pointed}\label{pro:nichols-twist}
There is a linear isomorphism
\[
\psi:\B(\hat{V})\to \B(V)
\]
such that $\psi(\hat{x}_i)= x_i$, $i\in \I_\theta$. 
Moreover, 
\begin{enumerate}[leftmargin=*]
\item[(i)]\cite[Lemma 4.4]{A-convex}
$\psi(\hat{x}_{\beta_M}^{n_M}\dots \hat{x}_{\beta_1}^{n_1})=f(\hat{x}_{\beta_M}^{n_M}\dots \hat{x}_{\beta_1}^{n_1})x_{\beta_M}^{n_M}\dots x_{\beta_1}^{n_1}$, for all $(n_i)_{i\in\I_M}\in \Z_{\geq 0}\a{M}$.
Thus,
\begin{align}\label{eqn:relation-products}
\psi(\hat{x}_\alpha^n)=\sigma(g_\alpha,g_\alpha)^{\binom{n}{2}}t_\alpha^n x_\alpha^n, \quad \alpha\in\Delta_+, n\in\N.
\end{align}
\item[(ii)]\cite[Lemma 1.1]{AS-annals} $\Delta(\psi(\hat{x}))=(\psi\ot\psi)\Delta_\sigma(\hat{x})$, $\hat{x}\in \B(\hat{V})$. 
In particular,
\begin{align}\label{eqn:relation-coproducts}
\Delta(x_\alpha^{N_\alpha})=\sigma(g_\alpha,g_\alpha)^{-\binom{N_\alpha}{2}}t_\alpha^{-N_\alpha}(\psi\ot\psi)\Delta_\sigma(\hat{x}_\alpha^{N_\alpha}), \quad \alpha\in\Delta_+.
\end{align}
\end{enumerate}
\qed
\end{pro}
We remark that we follow the notation of \cite[\S 4.1]{A-convex} rather than that of \cite{AS-pointed}, where the map $\B(\hat{V})\to \B(V)$ is denoted by $\varphi$.

\subsection{The strategy}\label{sec:strategy}

We give a brief overview of our strategy to compute all liftings of $V\in\ydh$ as cocycle deformations of $\B(V)\# H$.

We fix an {\it stratification} $\Gc=\Gc_0\sqcup\dots \sqcup \Gc_\ell$ of the set of generators $\Gc$ of the ideal defining $\B(V)$ and set $\B_0(V)=T(V)$, $\B_i(V)=T(V)/\lg \cup_{j=0}^{i-1}\Gc_j\rg$, $i>0$. In particular, $\Gc_i\subset \mP(\B_i)$, $i<\ell$, and $\k\lg \Gc_\ell\rg\subset \B_\ell(V)$ is a normal Hopf subalgebra. Also, we set $\mH_i=\B_i(V)\# H$, $i\geq 0$.

The program thus consists in $\ell+1$ steps. At each step $i\geq 0$, we have:
\begin{itemize} [leftmargin=*]\renewcommand{\labelitemi}{$\bullet$}
\item A collection $\Cleft'\mH_i$ of left $H$-module algebras $\mE_i(\bs\mu_{i-1})$, $\bs\mu_{i-1}\in\mR_{i-1}\subseteq \k^{\Gc_{i-1}}$, such that $\mA_i(\bs\mu_{i-1})\coloneqq \mE_i(\bs\mu_{i-1})\# H$ is a right cleft object of $\mH_i$, with coaction $\rho_i(\bs\mu_{i-1})$ and section $\gamma_i(\bs\mu_{i-1}):\mH_i\to \mA_i(\bs\mu_{i-1})$. 
\item A collection of cocycle deformations $\mL_i(\bs\mu_{i-1})$, each one equipped with a left coaction $\delta_i(\bs\mu_{i-1}):\mA_i(\bs\mu_{i-1})\to \mL_i(\bs\mu_{i-1})\ot \mA_i(\bs\mu_{i-1})$. 
\item A set of {\it deformation parameters}, see \cite[(3.9), (3.11)]{AAG}:
\begin{align}\label{eqn:R}
\mR_i=\{\bs\mu_i=(\mu_\alpha)_{\alpha\in\Gc_i}\in \k^{\Gc_i}|\mu_\alpha=0 \text{ if }  g_\alpha^{N_\alpha}=1 \text{ or } \chi_\alpha^{N_\alpha}\neq \eps\}.
 \end{align}
 \end{itemize}
Then step $i$ is an algorithmic computation that produces a new collection $\mE_{i+1}(\bs\mu_{i})$, $\mL_i(\bs\mu_{i})$, $\bs\mu_{i}\in\mR_{i}$, as above. We review this computation.

\begin{itemize} [leftmargin=*]\renewcommand{\labelitemi}{$\circ$}
\item We start with $\mR_{-1}=\{*\}$, $\mE_0=T(V)$, $\rho_0=\Delta$, $\gamma_0=\id$; hence $\mL_0=\mH_0$.

\item For each $\mA_i=\mE_i(\bs\mu_{i-1})\# H\in \Cleft'\mH_i$ and $\bs\mu_i\in\mR_i$, we set:
\begin{align}\label{eqn:def-E}
 \mE_{i+1}(\bs\mu_i) &=\mE_i/\lg \gamma_i(r)-\mu_r : r\in\Gc_i\rg
 \end{align} 
and we check that  $ \mE_{i+1}(\bs\mu_i)\neq 0$.

\item If  $ \mE_{i+1}(\bs\mu_i)\neq 0$, then we set 
\begin{align} \label{eqn:def-L} 
\mL_{i+1}(\bs\mu_i)&=\mL_i/\lg \tilde{r}-\mu_r(1-g_r): r\in\Gc_i\rg.
 \end{align} 
Here $\tilde{r}\in \mL_i$ is the unique solution of, \cf~\cite[(3.7)]{AAG}:
\begin{align}\label{eqn:solution1}
\tilde{r}\ot 1&=\delta_i(\gamma_i(r))-g_r\ot \gamma_i(r)\in \mL_i\ot \mA_i.
\end{align} 
Observe that $\gamma_i(r)$ can be computed recursively as the solution of:
\begin{align}\label{eqn:solution2}
\rho_i(\gamma_i(r))&=(\gamma_i\ot\id)\Delta(r)\in \mA_i\ot \mH_i.
\end{align}
 \end{itemize}
Hence step $i$ can be summarized as the determination of
 \begin{align*}
 \gamma_i(r), \text{ hence of }\rho_i(r), \quad \text{and} \quad  \delta_i(\gamma_i(r)), \qquad r\in\Gc_i,
 \end{align*} 
 together with $ \mE_{i+1}(\bs\mu_i)\neq 0$.
\begin{exa}\label{exa:step0}
For each $\bs\mu_0\in\mR_0$, we have
\begin{align*}
 \mE_{1}(\bs\mu_0) &=T(V)/\lg r-\mu_r : r\in\Gc_0\rg, \\
  \mL_{1}(\bs\mu_0)&=T(V)\# H/\lg r-\mu_r(1-g_r) : r\in\Gc_0\rg.
 \end{align*} 
\end{exa}
\begin{pro}\label{pro:strategy}\cite[Proposition 3.3]{AAG}
Fix $i\geq 0$, $\bs\mu_i\in\mR_i$, assume $\mE_{i+1}(\bs\mu_i)\neq 0$. Then
\begin{enumerate}
\item $\mA_{i+1}(\bs\mu_i)\in\Cleft\mH_{i+1}$.
\item $\mA_{i+1}(\bs\mu_i)$ is a ($\mL_{i+1}(\bs\mu_i)$,$\mH_{i+1}$) bi-cleft object.
\item $\mL_{i+1}(\bs\mu_i)$ is a cocycle deformation of $\mH_{i+1}$.
\end{enumerate}
If $i=\ell$, then $\gr \mL_{\ell+1}(\bs\mu_\ell)\simeq \B(V)\# H$.\qed
\end{pro}

We set $\mR\coloneqq \mR_0\times \dots\times \mR_\ell$, $\bs\mu=\bs\mu_0\times \dots\times \bs\mu_\ell\in\mR$. We assume that 
\begin{align}\label{eqn:assumption-neq-zero}
\mE_{i+1}(\bs\mu_i)\neq 0, \quad \forall\,i\geq 0.
\end{align}
We set $\u_{i+1}(\bs\mu)=\mL_{i+1}(\bs\mu_i)$, $i\geq 0$, and 
\[
\u(\bs\mu)\coloneqq\u_{\ell+1}(\bs\mu).
\]
This defines a family of liftings which are cocycle deformations of $\B(V)\# H$.

\begin{theorem}\label{thm:strategy}\cite[Theorem 3.5]{AAG}
Assume that \eqref{eqn:assumption-neq-zero} holds for all $\bs\mu\in\mR$. 
If $L$ is a lifting of $V\in\ydh$, then there is $\bs\mu\in\mR$ such that $L\simeq \u(\bs\mu)$.\qed
\end{theorem}

\begin{rem}\label{rem:strategy}
If the subalgebra $\k\lg\Gc_i\rg$ of $\mH_i$ generated by $\Gc_i$ is normal, then $\mE_{i+1}(\bs\mu_i)\neq 0$, for every $\bs\mu_i\in\mR_i$. See \cite[Theorem 3.1]{AAG}, \cite[Theorem 4]{Gu}.
\end{rem}

\subsubsection{Notation} 
We remark that the starting point of these computations is the fact that all of these objects $\mH_i,\mA_i,\mL_i$ are algebra quotients of $T(V)\# H$. We identify $H$ with a subalgebra, resp. Hopf subalgebra, of $\mA_i$, resp. $\mH_i$ and $\mL_i$. We denote the elements in these subalgebras by $h\in H$, indistinctly. 
 
In turn, we identify $V$ with a subspace in each of these algebras, namely the image of $V$ via the projection. We denote the (image of the) elements in $V$ by $x_j\in\mH_i$, $y_j\in\mA_i$ and $a_j\in\mL_i$, $j\in\I_\theta$. It follows that
 \[
 \rho_i(y_j)=y_j\ot 1+g_j\ot x_j, \quad \delta_i(y_j)=a_j\ot 1+g_j\ot y_j, \qquad j\in\I_\theta.
 \] 
 We write $x_\alpha,y_\alpha,a_\alpha$, $\alpha\in\Delta_+$, for the images of the root vectors on each quotient.
  
\section{Liftings of generic Cartan type}\label{sec:two-steps}

We review the results of the algorithm proposed in \cite{AAG} and described in \S \ref{sec:strategy} to compute the liftings of $V\in\ydh$ in a specially convenient setting, which is, however, quite ubiquitous, see Corollary \ref{cor:AS}. Recall that $\Gc$ stands for a set of generators of the ideal defining $\B(V)$.
We shall make the following two assumptions.

\begin{itemize} [leftmargin=*]
\item[({\it a})] There is a (2-step) stratification $\Gc=\Gc_0\sqcup \Gc_1$, with 
\begin{align}\label{eqn:2steps-1}
\Gc_0&=
\{\text{(generalized) quantum Serre relations}\}, 
& \Gc_1&=\mO(V).
\end{align}
\item[({\it b})] The (generalized) quantum Serre relations hold in any lifting, that is,
\begin{align}\label{eqn:2steps-2}
g_r=1 \text{ or } \chi_r\neq \eps, \ \forall\,r\in \Gc_0.
\end{align}
\end{itemize}
Thus, $\widetilde{\B}(V)=T(V)/\lg \Gc_0\rg$ by ({\it a}). We set $\mT(V)=T(V)\#H$, $\mH'=\widetilde{\B}(V)\# H$ and $\mH=\B(V)\#H$. As well, $\pi':\mT(V)\to \mH'$ is the canonical projection.

By hypothesis ({\it b}), step 0 is automatic: $\mR_0=\{(0,\dots,0)\}\subset \k^{\Gc_0}$. We set $\mE'=\widetilde{\B}(V)$ and let $\mA'=\mE'\# H$ be the trivial cleft object, \cf Example \ref{exa:step0}. All liftings are found after step 1, \cf \S \ref{sec:strategy}. 
We follow \eqref{eqn:R}, \eqref{eqn:def-E} and set
 \begin{align*}
 \mR&=\{\bs\mu=(\mu_\alpha)_{\alpha\in\mO(V)}\in \k^{\mO(V)}|\mu_\alpha=0 \text{ if }  g_\alpha^{N_\alpha}=1 \text{ or } \chi_\alpha^{N_\alpha}\neq \eps\}\\
\mE(\bs\mu)&:=\mE'/\lg \{ y_\alpha^{N_\alpha}-\mu_\alpha |\alpha\in\mO(V)\}\rg, \qquad  \mA(\bs\mu):=\mE(\bs\mu)\# H.
 \end{align*}
Then $\mA=\mA(\bs\mu)$ is a cleft object of $\mH$; we denote by $\gamma:\mH\to \mA$ the section and by $\rho:\mA \to \mA\ot\mH$ the coaction. Also, we denote by 
\begin{align*}
\tau_{\bs\mu}:\mA'\twoheadrightarrow \mA(\bs\mu)
\end{align*}
the natural algebra projection. Now, let $s(\alpha)\in\mH'$, $\alpha\in\Delta_+$, denote the unique solutions of:
\begin{align}\label{eqn:a-tilde}
s(\alpha)\ot 1&=(\pi'\ot \tau_{\bs\mu})\Delta(x_{\alpha}^{N_{\alpha}})
\end{align}
and set $\uvi_\alpha(\bs\mu):=s(\alpha)-a_\alpha^{N_\alpha}-\mu_\alpha g_\alpha^{N_\alpha}$, \cf \eqref{eqn:solution1}.

Next we follow \eqref{eqn:def-L} and let $\u(\bs\mu)$ be the quotient of $T(V)\#H$ by the ideal $I(\bs\mu)$ generated by the relations
	\begin{align}\label{eqn:presentation-generic}
r&=0, & r &\in\Gc_0;\\
\notag a_\alpha^{N_\alpha}&=\mu_\alpha(1-g_\alpha^{N_\alpha})-\uvi_\alpha(\bs\mu), & \alpha&\in\mO(V).
	\end{align}

\smallbreak
 	
In this setting, using Theorem \ref{thm:normal}, Theorem \ref{thm:strategy}, together with Remark \ref{rem:strategy}, implies the following.
\begin{theorem}\label{thm:2steps}
$I(\bs\mu)$ is a Hopf ideal and the Hopf algebra $\u(\bs\mu)$ is a lifting of $V$ for every $\bs\mu\in\mR$. Moreover, $\u(\bs\mu)$ is a cocycle deformation of $\B(V)\# H$.

Conversely, if  ({\it a}) and  ({\it b}) hold and $L$ is a lifting of $V$, then there is $\bs\mu\in\mR$ such that $L\simeq \u(\bs\mu)$.
In particular, every such lifting is a cocycle deformation of $\B(V)\# H$.\qed
\end{theorem}

\begin{rem}
Recall that $\uvi_\alpha(\bs\mu)\in\k\Gamma$, see (the proof of) \cite[Theorem 3.5]{AAG}. We shall give a description of these elements in Proposition \ref{pro:formulae}.
\end{rem}

\begin{cor}\label{cor:AS}
Let $V$ be a braided vector space of finite Cartan type.  Assume that $\Gamma$ as in \eqref{eqn:Gamma} is finite and all prime divisors of $|\Gamma|$ are $>7$. 

Let $L$ be a Hopf algebra whose infinitesimal braiding is a   principal realization $V\in\ydh$. Then there is $\bs\mu\in\k^{\Gc_1}$  such that $L\simeq \u(\bs\mu)$

In particular, $L$ is a cocycle deformation of $\B(V)\# H$.
\end{cor}
\pf
In this setting,  ({\it a}) and  ({\it b})  hold by \cite[Lemma 5.4]{AS-annals}.
\epf

\begin{rem}
If $H=\k\Gamma$, then this is \cite[Classification Theorem 0.1]{AS-annals}, using \cite[Theorem 2]{A-nichols}. This also extends \cite[Theorem 7.8]{Ma}.
\end{rem}

\subsection{On the coproduct of powers of Cartan root vectors}\label{sec:copro}
By \eqref{eqn:a-tilde} and \eqref{eqn:presentation-generic}, we see that in order to have an explicit description of a big family of liftings, it is necessary to have an explicit formula of the coproducts $\Delta(x_\alpha^{N_\alpha})$ of the powers of the root vectors. We investigate this in this section.

We need to fix some notation. Recall that $\Delta_+=\{\beta_1,\dots,\beta_M\}$ denotes the ordered set of positive roots, with 
$\beta_{i+1}>\beta_i$, $i\in\I_{M-1}$.

If $\alpha\in\mO=\Delta_+$, then we set $\bx_\alpha:=x_\alpha^{N_\alpha}$, $\bg_\alpha:=g_\alpha^{N_\alpha}$, $\bchi_\alpha:=\chi_\alpha^{N_\alpha}$.
We set 
\begin{align*}
\bx&:=\bx_{\beta_M}\dots \bx_{\beta_1}, & \bx^{\bf n}&:=\bx_{\beta_M}^{n_M}\dots \bx_{\beta_1}^{n_1},
\end{align*}
for ${\bf n}=(n_1,\dots,n_M)\in\Z^M$, $n_i\geq 0$, $i\in \I_M$. Also, $\bg, \bg^{\bf n}\in\Gamma$, $\bchi, \bchi_{\bbeta}^{\bf n}\in\Alg(H,\k)$ are set accordingly. 
 If ${\bf n}=(n_1,\dots,n_M)\in\Z^M$, we write
\[
\underline{\bf n}\coloneqq n_1\beta_1+\dots+ n_M\beta_M.
\]
By \cite[Theorem 4.13]{A-distinguished}, the subspace $Z(V)$ generated by the powers of root vectors is a braided sub-Hopf algebra; hence there are univocally determined scalars $r_{{\bf n},{\bf m}}(\alpha)\in\k$ such that
\begin{align}\label{eqn:coproduct}
\Delta(\bx_\alpha)=\bx_\alpha\ot 1+ 1\ot \bx_\alpha
+\sum_{\mathclap{\underline{\bf n}+\underline{\bf m}=N_\alpha \alpha }}r_{{\bf n},{\bf m}}(\alpha) \ \bx^{\bf n}\ot \bx^{\bf m}.
\end{align}
See \S \ref{sec:conj} for more details on the scalars $r_{{\bf n},{\bf m}}(\alpha)$.
\begin{pro}\label{pro:formulae}
For each $\alpha\in\mO(V)$, there are scalars $\mu({\bf m})\in\k$ such that
\begin{align}\label{eqn:uvi}
	\uvi_\alpha(\bs\mu)=\sum_{\mathclap{\underline{\bf n}+\underline{\bf m}=N_\alpha \alpha }}\mu_\alpha({\bf m})r_{{\bf n},{\bf m}}(\alpha)  
	a_{\beta_M}^{n_MN_{\beta_M}}\cdots a_{\beta_1}^{n_1N_{\beta_1}} g_{\beta_M}^{m_MN_{\beta_M}}\cdots g_{\beta_1}^{m_1N_{\beta_1}}.
\end{align}
\end{pro}
\pf
Follows by plugging  \eqref{eqn:coproduct} into \eqref{eqn:a-tilde}; so $\mu_\alpha({\bf m})=\tau_{\bs\mu}(\bx^{\bf m})$.
\epf

\subsection{On coproduct formulae}\label{sec:scalars}
We assume for the rest of this section that the braiding matrix $(\hat{q}_{ij})_{i,j\in\I}$ is such that
\begin{align}\label{eqn:hypothesis-AAR}
\hat{q}_{\alpha,\beta}^{N_\beta}=1, \quad \alpha,\beta\in\mO(V).
\end{align} 

In this setting, we provide an explicit general description of the scalars $r_{{\bf n},{\bf m}}$ involved in the formulae of the coproducts of the powers of the root vectors $x_\alpha^{N_\alpha}$, $\alpha\in\Delta_+$. Our computations rely on the ideas developed in \cite{AAR2}, where the matrix $(\hat{q}_{ij})_{i,j\in\I}$ is assumed to satisfy \eqref{eqn:hypothesis-AAR}.

\begin{rem}
We remark that \eqref{eqn:hypothesis-AAR} means no restriction to compute coproducts of root vectors, as any braiding $(q_{ij})_{i,j\in\I}$ is twist equivalent to a braiding $(\hat{q}_{ij})_{i,j\in\I}$ satisfying this condition. Hence, a coproduct formula for this case can be translated into one for the general case via Proposition \ref{pro:nichols-twist}.
\end{rem}

It follows from \eqref{eqn:hypothesis-AAR} that the braided Hopf algebra $\z(V)$ as in \eqref{eqn:sequence} becomes a usual Hopf algebra; thus it is isomorphic to the universal enveloping algebra $U(\n)$ of the Lie algebra 
\begin{align}\label{eqn:lie}
\n=\mP(\z(V)).
\end{align}
 Moreover, a basis of $\n$ is given by the images $\xi_\beta=\iota^*(y_\beta^{(N_\beta)})$, $\beta\in\mO(V)$, of the divided powers that generate $\L(V)$ \cf~\eqref{eqn:divided-powers}.  We recall that:
\begin{align}\label{eqn:potencia}
\xi_\beta^n=\iota^*(y_\beta^{(N_\beta)})^n=n!\iota^*(y_\beta^{(nN_\beta)}),
\end{align} 
that is, $\lg \xi_\beta^n,x_\beta^{nN_\beta}\rg=n!$; see \cite[Corollary 4.5]{AAR1} for details. 

Now, let $\Delta^{\n}$ be the root system associated to $\n$: each $\beta\in\Delta_+^{\n}$ arises as $N_\alpha\alpha$ for some $\alpha\in\Delta_+^{\q}$; thus $|\Delta_+^{\n}|=|\mO|$. A key step is the following: the roots in $\Delta_+^{\n}$, hence the set of generators $\{\xi_\beta\}_{\beta\in\mO}$ inherits an ordering from that of $\Delta_+$ in such a way that 
\begin{align}\label{eqn:order}
\xi_\beta < \xi_{\beta'}\quad \text{ if and only if } \quad  \beta'<\beta.
\end{align}

If $\beta,\gamma\in \mO(V)$, then it follows from \cite[Theorem 4.7]{AAR1} that there is a unique $\eta\in \Delta_+^{\n}$, associated to $N_\beta\beta+N_\gamma\gamma$, such that $[\xi_\beta,\xi_\gamma]\in\n_{\eta}$ and hence there is $c(\beta,\gamma)\in\k$ such that
\begin{align}\label{eqn:eigenspace}
[\xi_\beta,\xi_\gamma]=c(\beta,\gamma)\xi_{\eta},
\end{align}
as all of these eigenspaces are either zero or one-dimensional.

Now, to determine the scalars $(r_{{\bf n},{\bf m}})_{{\bf n},{\bf m}}\in\k$ in \eqref{eqn:coproduct}, we proceed as follows. 
On the one hand, using \eqref{eqn:potencia}, \cf also \eqref{eqn:convolution}:
\begin{align}\label{eqn:evaluation}
\prod_{i=1}^{k}n_i!\prod_{j=1}^{l}m_j!\,r_{{\bf n},{\bf m}}&=\left\lg \xi_{\gamma_l}^{m_l}\dots \xi_{\gamma_1}\a{m_1} \xi_{\beta_k}\a{n_k}\dots\xi_{\beta_1}\a{n_1} ,x_\alpha^{N_\alpha}\right\rg.
\end{align}
Next we write the element $\xi_{\gamma_l}\a{m_l}\dots \xi_{\beta_1}^{n_1}$ in the PBW basis of $U(\n)$; observe that $\xi_{\gamma_j}<\xi_{\beta_i}$, for every $i,j$ by \eqref{eqn:order}. Now, an iterative application of \eqref{eqn:eigenspace} defines a scalar $c(\alpha)\in\k$ such that
\begin{align}\label{eqn:c-omega}
\xi_{\gamma_l}\a{m_l}\dots \xi_{\beta_1}^{n_1}=c(\alpha)\xi_\alpha + \text{other monomials}.
\end{align}
\begin{lem}\label{lem:zero}
$\left\lg \xi_{\gamma_l}^{m_l}\dots \xi_{\gamma_1}\a{m_1} \xi_{\beta_k}\a{n_k}\dots\xi_{\beta_1}\a{n_1} ,x_\alpha^{N_\alpha}\right\rg=c(\alpha)$. 
\end{lem}
\pf
Set $T_\alpha=\{\bs\tau=(\tau_1,\dots,\tau_s)\in\mO^s|\tau_1<\dots<\tau_s, s\in \N \}$ and let $(f_{\bs\tau,{\bs t}})_{\bs\tau,{\bf t}}\in\k$ be such that
\begin{align*}
\xi_{\gamma_l}\a{m_l}\dots \xi_{\beta_1}^{n_1}=c(\alpha)\xi_\alpha + 
\sum_{\mathclap{\bs\tau\in T_\alpha,{\bs t}\in\N^{|\bs\tau|}}}\ f_{\bs\tau,{\bs t}}\xi_{\tau_1}^{t_1}\dots \xi_{\tau_s}\a{t_s},
\end{align*}
\cf \eqref{eqn:c-omega}. In particular, $f_{\bs\tau,{\bs t}}\neq 0$ only if $\sum_i t_i\tau_i=N_\alpha\alpha$. We claim that 
\begin{align}\label{eqn:zero}
\lg\xi_{\tau_1}^{t_1}\dots \xi_{\tau_s}\a{t_s},x_\alpha\a{N_\alpha}\rg=0.
\end{align}
for every $\bs\tau\in T_\alpha$, ${\bf t}\in\N\a{\bs\tau}$. This implies the lemma, as $\lg\xi_\alpha,x_\alpha\a{N_\alpha}\rg=1$.

Let $\Delta^{(s)}$, $s\geq 1$, denote the $s$-th iteration of $\Delta=\Delta\a{(1)}$. We observe that \eqref{eqn:zero} is nonzero only if a summand of the form 
$x_{\tau_s}\a{r_sN_{\tau_s}} \ot \cdots \ot x_{\tau_1}\a{r_1N_{\tau_1}}$ appears when computing $\Delta^{(s-1)}(x_\alpha\a{N_\alpha})$, but the existence of such a term contradicts (an iteration of) Lemma \ref{lem:coproduct}.
\epf
We arrive to the following.
\begin{pro}\label{pro:scalars}
Assume that \eqref{eqn:hypothesis-AAR} holds. If $r_{{\bf n},{\bf m}}(\alpha)$ is as in \eqref{eqn:coproduct}, then
\begin{align}\label{eqn:c-t}
r_{{\bf n},{\bf m}}(\alpha)&=\prod_{i=1}^{k}\frac{1}{n_i!}\prod_{j=1}^{l}\frac{1}{m_j!}\,c(\alpha).
\end{align}
\end{pro}
\pf
Combine \eqref{eqn:evaluation} with Lemma \ref{lem:zero}.
\epf
Hence, we can determine all scalars $(r_{{\bf n},{\bf m}})$ out of the Lie algebra structure of $\n$; see the proof of Proposition \ref{pro:scalarsG2} for a concrete example.

\subsection{A conjecture}\label{sec:conj}

It follows from \cite[Theorem 4.13]{A-distinguished} that 
\[
t_{{\bf n},{\bf m}}=0 \text{ or } n_i=0 \text{ for each }\beta_i>\alpha.
\]
In plain words, only root vectors associated with roots lesser than $\alpha$ may appear in the left tensorand. We believe there is a mirror situation for ${\bf m}$, in the sense that only with roots bigger than $\alpha$ may appear; see Remark \ref{rem:conj}.
That is, we conjecture the following: set 
\begin{align}\label{eqn:omega-alpha}
\Omega_{\alpha}&=\{({\bf n},{\bf m})\in \Z^{2M} | n_i=0\text{ if }\beta_i>\alpha, m_i=0\text{ if }\beta_i<\alpha; i\in\I_M\}.
\end{align}

\begin{conjecture}\label{conj-1}
There are scalars $(r_{{\bf n},{\bf m}})_{{(\bf n,m)}\in\Omega_{\alpha}}\in\k$ such that
\begin{align*}
\Delta(\bx_\alpha)=\bx_\alpha\ot 1+ 1\ot \bx_\alpha
+\sum_{\mathclap{\underline{\bf n}+\underline{\bf m}=N_\alpha \alpha }}r_{{\bf n},{\bf m}} \ \bx^{\bf n}\ot \bx^{\bf m}.
\end{align*}
\end{conjecture}
Moreover, we further propose the following. We write 
\begin{align*}
\underline{x}&=x_{\beta_M}\dots x_{\beta_1}, & \underline{x}^{\bf n}&=x_{\beta_M}^{n_m}\dots x_{\beta_1}^{n_1}, \quad {\bf n}=(n_i)_{i\in\I_M}\in\Z^M.
 \end{align*}

\begin{conjecture}\label{conj-2}
There are scalars $(s_{{\bf n},{\bf m}})_{{(\bf n,m)}\in\Omega_{\alpha}}\in\k$ such that
\begin{align*}
\Delta(x_\alpha)&=x_\alpha\ot 1+ 1\ot x_\alpha+\sum_{\mathclap{\underline{\bf n}+\underline{\bf m}= \alpha }} \ s_{{\bf n},{\bf m}}\,
 \underline{x}^{\bf n}  \ot \underline{x}^{\bf m}.
 \end{align*}
\end{conjecture}

The following lemma is now a direct consequence of \cite[Theorem 4.13]{A-distinguished}.
\begin{lem}\label{lem:coproduct}
If Conjecture \ref{conj-2} holds, then Conjecture \ref{conj-1} holds. 
\end{lem}
\pf
We recall from \cite[Lemma 4.5]{A-convex} the braided commutator between two root vectors is a linear combination of products of root vectors associated to intermediate roots.
Hence the resulting basic tensors in the computation of $\Delta(x_\alpha)^{N_\alpha}$ are such that left tensorands are products of root vectors associated to roots lesser or equal than $\alpha$, while the right ones are products of root vectors associated to roots bigger or equal than $\alpha$. 

The lemma now follows by \eqref{eqn:coproduct}, as the only possibilities for a basic tensor to involve the root $\alpha$ are $x_\alpha^{N_\alpha}\ot 1$ and $1\ot x_\alpha^{N_\alpha}$.
\epf

\begin{rem}\label{rem:conj}
Conjecture \ref{conj-2}, and hence \ref{conj-1}, holds when 
\begin{enumerate}
\item All root vectors are primitive in $\widetilde{\B}(V)$.
\item $V$ is of Cartan type $A_n$, $n\in\N$. See \cite[Lemma 6.9]{AS-pointed}, also \cite{AD,AK}.
\item $V$ is of Cartan type $B_n$, $C_n$, $D_n$, $n\in\N$. See \cite{K1,K2}.
\item $V$ is of Cartan type $G_2$. See \cite{KV}, also Proposition \ref{pro:copro-g2}.
\item $V$ is of super type $A_n$, $n\in\N$. See \cite[Lemma 34]{A-distinguished}.
\item $V$ is of modular type $\mathfrak{br}(2;5)$. See the proof of \cite[Lemma 37]{A-distinguished}.
\item $V$ is of unidentified type $\mathfrak{ufo}(7)$.
\end{enumerate}
Case (7) is a particular case of (1). For (3), the coproduct formula is not in terms of the root vectors, but rather of a linear combination of them, suitable for the purposes in the references. Idem for case (4) in \cite{KV}; hence the need of Proposition \ref{pro:copro-g2}.
\end{rem}

\section{Cartan Type  $G_{2}$}\label{sec:cartan}

We fix a braided vector space $V$ of diagonal type Cartan $G_2$ with parameter $q\in\G_N'$, $N>3$; see \eqref{eqn:braiding-G2}. The braiding of $V$ is encoded in the {\it generalized Dynkin diagram} associated to $V$:
\begin{align}\label{eq:dynkin-type-G}
\xymatrix{  \overset{q}{\underset{\ }{\circ}} \ar  @{-}[r]^{q^{-3}} &
\overset{\,\,q^3}{\underset{\ }{\circ}}}
\end{align}
We fix a cosemisimple Hopf algebra $H$ with a principal realization $V\in\ydh$ and $\Gamma\leq Z(G(H))$ as in \eqref{eqn:Gamma}. 

The Nichols algebra $\B(V)$ is generated by $x_1, x_2$ with defining
relations
\begin{align}\label{eq:rels-type-G-1}
[x_{12},x_2]_c &= 0; \qquad x_{11112}=0;\\
\label{eq:rels-type-G-2} x_{\alpha}^{N_\alpha}&=0, \qquad 
\alpha\in\Delta_+.
\end{align}
Recall that the relations in \eqref{eq:rels-type-G-1} are the {\it quantum Serre relations} while the ones in \eqref{eq:rels-type-G-2} are the {\it powers of root vectors}, indexed by the (positive) roots $$\Delta_+=\{\alpha_1,\alpha_2,\alpha_1+\alpha_2,2\alpha_1+\alpha_2,3\alpha_1+\alpha_2,3\alpha_1+2\alpha_2\}.$$
In this case, we have $\Delta_+=\mO(V)$, that is all roots are Cartan.
More precisely, the root vectors are
\begin{align*}
x_{\alpha_1}&=x_1, & x_{\alpha_2}&=x_2, & x_{\alpha_1+\alpha_2}&=x_{12},\\
x_{2\alpha_1+\alpha_2}&=x_{112}, & x_{3\alpha_1+\alpha_2}&=x_{1112} & x_{3\alpha_1+2\alpha_2}&=[x_{112},x_{12}]_c.
\end{align*}
We shall denote $\beta:=3\alpha_1+2\alpha_2$ for shortness. We will also refer indistinctly to the roots by their indices, that is $\Delta_+=\{1,2,12,112,1112,11212\}$. We recall that there is a preferred order of the roots, hence of the root vectors, namely the {\it convex order}
\begin{align}\label{eqn:order-G}
\alpha_1<3\alpha_1+\alpha_2<2\alpha_1+\alpha_2<3\alpha_1+2\alpha_2<\alpha_1+\alpha_2<\alpha_2.
\end{align}

Set $M=\dfrac{N}{(N,3)}$. We have, cf.~\eqref{eq:rels-type-G-2}:
\begin{align}\label{eqn:N_alpha}
\begin{split}
& N_{\alpha_1}=N_{\alpha_1+\alpha_2}=N_{2\alpha_1+\alpha_2}=N, \\
& N_{\alpha_2}=N_{3\alpha_1+\alpha_2}=N_{3\alpha_1+2\alpha_2}=M.
\end{split}
\end{align}
We fix an stratification $\Gc=\Gc_0\sqcup\Gc_1$ as
\begin{align*}
\Gc_0&=\{\text{quantum Serre relations }\eqref{eq:rels-type-G-1}\},\\
\Gc_1&=\{\text{powers of root vectors }\eqref{eq:rels-type-G-2}\}.
\end{align*}
We set $\mH'=\mT(V)/\lg\Gc_0\rg=\widetilde{\B}(V)\#H$, $\mH=\B(V)\#H=\mH'/\lg\Gc_1\rg$.

We remark that there are two distinguished cases, according to whether the quantum Serre relations \eqref{eq:rels-type-G-1} can be deformed or else hold in any lifting. We have, {\it cf.} \cite[Proposition 3.2]{AGI} also \cite[Theorem 5.6]{AS-pointed}, that these relations can be deformed 
if and only if:
\begin{align}\label{eqn:quantum-def}
\qb=\qb_d:=\begin{pmatrix}
q&q^3\\ q&q^3
\end{pmatrix} \quad \text{(hence $N=7$)}.
\end{align}
Hence, we shall consider separately
\begin{itemize}
\item {\it The generic case:} when $\qb\neq \qb_d$. See \S \ref{sec:generic}.
\item {\it The degenerate case:} when $\qb= \qb_d$. See \S \ref{sec:deformed}.
\end{itemize}

\subsection{The generic case}\label{sec:generic}
This case agrees with the setting of \S \ref{sec:two-steps}, hence all liftings arise as $\u(\bs\lambda)$ in Theorem \ref{thm:2steps}. To describe them, we need an explicit description of the coproducts of the powers of the root vectors. For this we use Proposition \ref{pro:nichols-twist} and the results in \cite{AAR2}.

We start with a lemma in which we exhibit the explicit formulae which describe the brackets between two root vectors.
 \begin{lem}\label{lem:convex}
On the one hand,
\begin{align}\label{eqn:bracket-2}
\begin{split}
[x_1,x_2]_c&=x_{12};  \qquad \quad \ \, [x_{1},x_{12}]_c=x_{112}; \\
 [x_{1},x_{112}]_c&=x_{1112}; \qquad  [x_{112},x_{12}]_c=x_\beta.
\end{split}
\end{align}
Also, there is a subset of $q$-commuting vectors:
\begin{align}\label{eqn:bracket-1}
[x_{12},x_2]_c&=0; &[x_{\beta},x_{12}]_c&=0; &[x_{112}, x_{\beta}]_c&=0; \\
\notag  [x_{1112},x_{112}]_c&=0; &[x_{1},x_{1112}]_c&=0. &&
\end{align}
For the rest of the roots, we have:
\begin{align}\label{eqn:bracket-3}
[x_\beta,x_2]_c&=q_{12}^2q^3(q^2-1)(q-1)x_{12}^3; & [x_{112},x_2]_c&=q_{12}q(q^2-1)x_{12}^2;\\
\notag [x_{1112},x_{12}]_c&=\dfrac{q_{12}q(q^3-1)}{q+1}x_{112}^2; &  [x_{1},x_{\beta}]_c&=\dfrac{q_{12}q(q^3-1)}{q+1}x_{112}^2;\\
\notag [x_{1112},x_{\beta}]_c&=\dfrac{q^3q_{12}^2(q-1)(q^3-1)}{q+1}x_{112}^3; & &
\end{align}
\vspace*{-.4cm}
\begin{align*}
\notag [x_{1112},x_{2}]_c&=
q_{12}q(q^2-q-1)x_{\beta} +q_{12}^2q^2(q^3-1)x_{12}x_{112}. \ \ \ \ \ \ \ \ \ \ \ \ \ \ \
\end{align*}
\end{lem}
\pf
First, \eqref{eqn:bracket-2} is by definition. Now, we recall from \cite[Lemma 4.5]{A-convex} that  the bracket of two given root vectors is a linear combination of intermediate root vectors. In particular, 
the bracket of consecutive root vectors is zero, this shows \eqref{eqn:bracket-1}, \cf~\eqref{eqn:order-G}.  Now 
\eqref{eqn:bracket-3} follows by an easy computation.
\epf

We continue with a well-known result, \cf~\cite[Theorem 4.2]{KV}.
\begin{pro}\label{pro:copro-g2}
In $\mH'$, we have:
\begin{align*}
\Delta(x_i)&=x_i\ot 1+g_i\ot x_i, \quad i=1,2;\\
\Delta(x_{12})&=x_{12}\ot 1+g_{12}\ot x_{12}+(1-q^{-3})x_1g_2\ot x_2;\\
\Delta(x_{\beta})&=x_{\beta}\ot 1+g_{\beta}\ot x_{\beta}+ q^2(1-q^{-3})^2x_{112}x_1g_2\ot x_2 \\
&\quad +q_{21}(1-q^{-3})(q^3-q^2-q)x_{1112}g_2\ot x_2 \\
&\quad +q_{21}(1-q^{-3})^2(1-q^{-2})(1-q^{-1})x_1^3g_{22}\ot x_2^2 \\
&\quad +q^2(1-q^{-3})^2(1-q^{-2})x_1^2g_{122}\ot x_2x_{12} \\
&\quad +q^2(1-q^{-3})x_{112}g_{12}\ot x_{12} +q^2(1-q^{-3})(1-q^{-2})x_1g_{1122}\ot x_{12}^2; \\
\Delta(x_{112})&=x_{112}\ot 1+g_{112}\ot x_{112}+(1-q^{-3})(1-q^{-2})x_1^2g_2\ot x_2 \\
&\quad +(1-q^{-2})(1+q)x_1g_{12}\ot x_{12};\\
\Delta(x_{1112})&=x_{1112}\ot 1+g_{1112}\ot x_{1112} \\
&\quad +(1-q^{-3})(1-q^{-2})(1-q^{-1})x_1^3g_2\ot x_2 \\
&\quad +(1-q^{-3})(1-q^{-2})q^2x_1^2g_{12}\ot x_{12} + (1-q^{-3})q^2x_1g_{112}\ot x_{112}.
\end{align*}
\end{pro}
\pf
This is a straightforward computation, using Lemma \ref{lem:convex}.
\epf

In particular, Proposition \ref{pro:copro-g2} shows that Conjecture \ref{conj-2} holds and therefore there are scalars $(r_{{\bf n},{\bf m}})_{{({\bf n},{\bf m})}\in\Omega_{\alpha}}$ such that the coproducts of the powers of the root vectors follow the rule in Conjecture \ref{conj-1}.

 For simplicity, and following the notation in \cite[\S 4.]{AAR2}, we shall denote these scalars $r_{{\bf n},{\bf m}}$ by $(\ta_i)_{i\in\I_{12}}$ (when $3\nmid N$) and $(\tb_i)_{i\in\I_{12}}$ (when $3|N$). Hence

$\bullet$ When $ (N, 3) = 1 $:\label{page:ai}
\begin{align*}
\Delta(x_{1}^{N})&=x_{1}^{N}\ot 1+g_{1}^{N}\ot x_{1}^{N}; \qquad \Delta(x_{2}^{N})=x_{2}^{N}\ot 1+g_{2}^{N}\ot x_{2}^{N};\\
\Delta(x_{12}^N) &= x_{12}^N\ot 1+g_{12}^N\ot x_{12}^N +\ta_1 x_1^Ng_2^N\ot x_2^N;\\
\Delta(x_{112}^N) &= x_{112}^N\ot 1+g_{112}^N\ot x_{112}^N  + \ta_2x_1^Ng_{12}^N\ot x_{12}^N + \ta_3x_1^{2N}g_2^N\ot x_2^N; \\
\Delta(x_{1112}^{N}) &= x_{1112}^{N}\ot 1+g_{1112}^{N}\ot x_{1112}^{N} + \ta_4x_1^Ng_{112}^{N}\ot x_{112}^{N} + \ta_5x_1^{2N}g_{12}^{N}\ot x_{12}^{N}\\
&\quad +\ta_6x_1^{3N}g_{2}^{N}\ot x_{2}^{N};\\
\Delta(x_{\beta}^{N}) &= x_{\beta}^{N}\ot 1+g_{\beta}^{N}\ot x_{\beta}^{N} + \ta_7x_{112}^Ng_{12}^N\ot x_{12}^N + \ta_8x_{1112}^{N}g_2^{N}\ot x_2^{N} \\
&\quad + \ta_9x_1^{N}g_{12}^{2N}\ot x_{12}^{2N} + \ta_{10}x_1^{2N}g_{122}^{N}\ot x_{2}^{N}x_{12}^N  \\
&\quad + \ta_{11}x_{112}^Nx_1^Ng_2^{N}\ot x_2^{N}+ \ta_{12}x_1^{3N}g_{2}^{2N}\ot x_{2}^{2N}. 
\end{align*}

$\bullet$ When $ (N, 3) = 3 $ (then $ N=3M $):
\begin{align*}
\Delta(x_{1}^{N})&=x_{1}^{N}\ot 1+g_{1}^{N}\ot x_{1}^{N}; \qquad \Delta(x_{2}^{M})=x_{2}^{M}\ot 1+g_{2}^{M}\ot x_{2}^{M};\\
\Delta(x_{12}^N)&=x_{12}^N\ot 1+g_{12}^N\ot x_{12}^N + \tb_1 x_{\beta}^{M}g_2^{M}\ot x_2^{M} + \tb_2x_{1112}^{M}g_2^{2M}\ot x_2^{2M} \\
&\quad +\tb_3x_1^Ng_2^N\ot x_2^N;\\
\Delta(x_{112}^N)&=x_{112}^N\ot 1+g_{112}^N\ot x_{112}^N + \tb_4x_1^Ng_{12}^N\ot x_{12}^N + \tb_5x_{1112}^{M}g_{\beta}^{M}\ot x_{\beta}^{M} \\
&\quad  + \tb_6x_1^{2N}g_2^N\ot x_2^N + \tb_7x_{1112}^{2M}g_{2}^{M}\ot x_{2}^{M} + \tb_8x_{1112}^{M}x_1^Ng_{2}^{2M}\ot x_{2}^{2M} \\
& \quad  + \tb_9x_1^Ng_{12}^N\ot x_{2}^{M}x_{\beta}^{M}; \\
\Delta(x_{1112}^{M})&=x_{1112}^{M}\ot 1+g_{1112}^{M}\ot x_{1112}^{M}+\tb_{10}x_1^Ng_2^{M}\ot x_2^{M}; \\
\Delta(x_{\beta}^{M})&=x_{\beta}^{M}\ot 1+g_{\beta}^{M}\ot x_{\beta}^{M} + \tb_{11}x_1^Ng_2^{2M}\ot x_2^{2M}
+\tb_{12}x_{1112}^{M}g_2^{M}\ot x_2^{M}. 
\end{align*}
We follow the ideas in \S \ref{sec:twist-equiv} and \S \ref{sec:scalars} to compute the scalars $\ta_i$ and $ \tb_i $.
\begin{pro}\label{pro:scalarsG2}
For the case $(N,3)=1$:
\begin{align*}
\ta_1&=(1-q^{-3})^Nq_{21}^{\frac{N(N-1)}{2}};\\
\ta_2&=2(1-q^{-2})^Nq_{21}^{\frac{N(N-1)}{2}};\\
\ta_3&=(1-q^{-3})^N(1-q^{-2})^Nq_{21}^{N(N-1)};\\
\ta_4&=3(1-q^{-1})^Nq_{21}^{\frac{N(N-1)}{2}};\\
\ta_5&=3(1-q^{-2})^N(1-q^{-1})^Nq_{21}^{N(N-1)};\\
\ta_6&=(1-q^{-3})^{N}(1-q^{-2})^{N}(1-q^{-1})^{N}q_{21}^{\frac{3N(N-1)}{2}};\\
\ta_7&=3(1-q^{-1})^{N}q_{21}^{\frac{N(N-1)}{2}};\\
\ta_8&=-(1-q^{-3})^{N}q_{21}^{\frac{N(3N-1)}{2}};\\
\ta_9&=3(1-q^{-2})^{N}(1-q^{-1})^{N}q_{21}^{N(N-1)};\\
\ta_{10}&=3(1-q^{-3})^{N}(1-q^{-2})^{N}(1-q^{-1})^{N}q_{21}^{\frac{3N(N-1)}{2}};\\
\ta_{11}&=3(1-q^{-3})^{N}(1-q^{-1})^{N}q_{21}^{N(N-1)};\\
\ta_{12}&=(1-q^{-3})^{2N}(1-q^{-2})^{N}(1-q^{-1})^{N}q_{21}^{N(3N-2)}.
\end{align*}
On the other hand, when $(N,3)=3$, we have:
\begin{align*}
\tb_1&=3(1-q^{-3})^M(1-q^{-2})^{-M}(1-q^{-1})^{-M}q_{21}^{\frac{M(N+1)}{2}};\\
\tb_2&=3(1-q^{-3})^{2M}(1-q^{-2})^{-M}(1-q^{-1})^{-M}q_{21}^{NM};\\
\tb_3&=(1-q^{-3})^Nq_{21}^{\frac{N(N-1)}{2}};\\
\tb_4&=-(1-q^{-2})^{N}(1-q^{-1})^{N}q_{21}^{\frac{N(N-1)}{2}};\\
\tb_5&=3(1-q^{-2})^{-M}(1-q^{-1})^{M}q_{21}^{\frac{M(N+1)}{2}};\\
\tb_6&=(1-q^{-3})^N(1-q^{-2})^Nq_{21}^{N(N-1)};\\
\tb_7&=3(1-q^{-3})^M(1-q^{-2})^{M}(1-q^{-1})^{M}q_{21}^{NM};\\
\tb_8&=3(1-q^{-3})^{2M}(1-q^{-2})^{2M}(1-q^{-1})^{2M}q_{21}^{\frac{N(N-1)}{2}};\\
\tb_9&=3(1-q^{-3})^M(1-q^{-2})^{2M}(1-q^{-1})^{2M}q_{21}^{M(N-1)};\\
\tb_{10}&=(1-q^{-3})^{M}(1-q^{-2})^{M}(1-q^{-1})^{M}q_{21}^{\frac{N(M-1)}{2}}; \\
\tb_{11}&=(1-q^{-3})^{2M}(1-q^{-2})^{M}(1-q^{-1})^{M}q_{21}^{M(N-2)};\\
\tb_{12}&=2(1-q^{-3})^{M}q_{21}^{\frac{M(N-1)}{2}}.
\end{align*}
\end{pro}
\pf 
First, we compute the scalar $\ta_1$. By Proposition \ref{pro:copro-g2}, 
\[
\Delta(x_{12}) = x_{12}\ot 1+g_{12}\ot x_{12}+(1-q^{-3})x_1g_2\ot x_2.\]
Hence, 
\begin{align*}
\Delta(x_{12}^N) & = \sum_{d_i\in \mathcal{D}_{12}} d_1 d_2 \cdots d_N, 
\end{align*} 
where $ \mathcal{D}_{12} = \{ x_{12}\ot 1, g_{12}\ot x_{12}, (1-q^{-3})x_1g_2\ot x_2 \} $. Using that the order is convex, the only way to have $ x_{2}^N $ in the second tensorand is by taking $ d_1 = \cdots = d_N = (1-q^{-3})x_1g_2\ot x_2 $. Therefore, the part that we focus on is:
\begin{align*}
\left((1-q^{-3})x_1g_2\ot x_2\right)^N &= (1-q^{-3})^N(x_1g_2)^N\ot x_2^N \\
&= (1-q^{-3})^Nq_{21}^{\frac{N(N - 1)}{2}}x_1^Ng_2^N\ot x_2^N
\end{align*}
and thus $\ta_1=(1-q^{-3})^Nq_{21}^{\frac{N(N-1)}{2}}$.
We observe that the scalars $\ta_{3}, \ta_{6}, \ta_{12}, \tb_{3}$, $\tb_{6}$, $\tb_{10} $ and $ \tb_{11} $ can be calculated in the same way.

For the remaining scalars, we use the techniques from \S\ref{sec:scalars}: we calculate them when \eqref{eqn:hypothesis-AAR} holds, that is when $q_{\alpha\alpha'}^{N_{\alpha'}}=1$, for all $\alpha,\alpha'$. Afterwards we use the formulae from \S \ref{sec:twist-equiv} to get the general case. We shall compute the scalar $\hat{\ta}_2$, the others follow similarly. 

We fix a braiding matrix $(\hat{q}_{ij})_{i,j\in\I}$ satisfying \eqref{eqn:hypothesis-AAR}. In this case, the Lie algebra $\n=\k\{\xi_1,\xi_2\}$ \cf \eqref{eqn:lie} is of type $G_2$, see \cite[\S 4 (Row 10)]{AAR2}. 

We fix $\xi_1>\xi_2$, \cf \eqref{eqn:order}. Now, we know that $ \xi_{12} = c_1[\xi_2, \xi_1] $ and $ \xi_{112} = c_2[\xi_{12}, \xi_{1}] $ for some $ c_i\in\Bbbk $. Therefore, in $ U(\n) $, we have:
\begin{align*}
\xi_{2}\xi_{1} &= c_1^{-1}\xi_{12}+\xi_{1}\xi_{2}; \\
\xi_{2}\xi_{1}^2 &= (c_1^{-1}\xi_{12}+\xi_{1}\xi_{2})\xi_{1} = c_1^{-1}(c_2^{-1}\xi_{112}+\xi_{1}\xi_{12})+c_1^{-1}\xi_{1}\xi_{12}+\xi_{1}^2\xi_{2} \\
&= c_1^{-1}c_2^{-1}\xi_{112}+2c_1^{-1}\xi_{1}\xi_{12}+\xi_{1}^2\xi_{2}; \\
\xi_{12}\xi_{1} &= c_2^{-1}\xi_{112}+\xi_{1}\xi_{12}.
\end{align*}
By \eqref{eqn:c-t}, we have that:
\begin{align*}
\hat{\ta}_1 &= c_1^{-1}; & \hat{\ta}_3 &= \frac{c_1^{-1}c_2^{-1}}{2!}; & \hat{\ta}_2 &= c_2^{-1}. 
\end{align*} 
As we have already computed $ \hat{\ta}_1 $ and $ \hat{\ta}_3 $, we obtain 
\[
\hat{\ta}_2=c_2^{-1}=2\hat{\ta}_3c_1=2\hat{\ta}_3\hat{\ta}_1^{-1}=2(1-q^{-2})^N\hat{q}_{21}^{\frac{N(N-1)}{2}}.
\]
To find the scalar $ \ta_2 $ in the general case, we follow the recipe from \S \ref{sec:twist-equiv}. We refer the reader to \loc for unexplained notation. To begin with, we need 
$\Delta_\sigma(\hat{x}_{112}^N) $. For this, we first observe that:
\begin{align*}
\Delta(\hat{x}_{112}^N) &= \hat{x}_{112}^N\ot 1 + g_{112}^N\ot \hat{x}_{112}^N + \hat{\ta}_2\hat{x}_1^Ng_{12}^N\ot \hat{x}_{12}^N + \hat{\ta}_3\hat{x}_1^{2N}g_{2}^{N}\ot \hat{x}_{2}^{N} 
\end{align*}
Then, by \eqref{eqn:coproddeformed}, we have
\begin{align*}
\Delta_\sigma(\hat{x}_{112}^N) &= \sigma(g_{112}^N, 1)\hat{x}_{112}^N\ot 1 + \sigma(1, g_{112}^N)g_{112}^N\ot \hat{x}_{112}^N \\
&\quad + \sigma(g_{1}^N, g_{12}^N)\hat{\ta}_2\hat{x}_1^Ng_{12}^N\ot \hat{x}_{12}^N + \sigma(g_{1}^{2N}, g_{2}^N)\hat{\ta}_3\hat{x}_1^{2N}g_{2}^{N}\ot \hat{x}_{2}^{N}\\
&= \hat{x}_{112}^N\ot 1 + g_{112}^N\ot \hat{x}_{112}^N + \sigma(g_{1}, g_{2})^{N^2}\hat{\ta}_2\hat{x}_1^Ng_{12}^N\ot \hat{x}_{12}^N \\
&\quad + \sigma(g_{1}, g_{2})^{2N^2}\hat{\ta}_3\hat{x}_1^{2N}g_{2}^{N}\ot \hat{x}_{2}^{N}.
\end{align*} 
Using Proposition \ref{pro:nichols-twist}~(i), we get that:
\begin{align*}
\psi(\hat{x}_{112}^N) &= \sigma(g_{1}, g_{2})^{N(N+1)} x_{112}^N; & \psi(\hat{x}_{12}^N) &= \sigma(g_{1}, g_{2})^{\frac{N(N+1)}{2}} x_{12}^N; \\
\psi(\hat{x}_{1}^N) &= x_{1}^N; & \psi(\hat{x}_{2}^N) &= x_{2}^N.
\end{align*}
We also have $ t_{112} = \sigma(g_{1}, g_{2})^2$, \cf~\eqref{eqn:t-alpha}. Finally, from \eqref{eqn:relation-coproducts}, we obtain: 
\begin{align*}
\Delta(x_{112}^N) &= \sigma(g_{112}, g_{112})^{-\frac{N(N-1)}{2}}t_{112}^{-N}(\psi\ot\psi)\Delta_\sigma(\hat{x}_{112}^N) \\
&=x_{112}^N\ot 1 + g_{112}^N\ot x_{112}^N + 2(1-q^{-2})^N q_{21}^{\frac{N(N-1)}{2}}x_1^Ng_{12}^N\ot x_{12}^N \\
&\quad + \ta_3x_1^{2N}g_{2}^{N}\ot x_{2}^{N}
\end{align*} 
and hence $\ta_2=2(1-q^{-2})^Nq_{21}^{\frac{N(N-1)}{2}}$.
\epf

Now, consider the set of deformation parameters, \cf\eqref{eqn:R}:
\begin{align}\label{eqn:mu}
\mR=\{\bs\mu=(\mu_\alpha)_{\alpha\in\Delta_+}\in \k^{\Delta_+}|\mu_\alpha=0 \text{ if }  g_\alpha^{N_\alpha}=1 \text{ or } \chi_\alpha^{N_\alpha}\neq \eps\}.
\end{align}

We follow \eqref{eqn:presentation-generic}, see also \eqref{eqn:uvi}, to give an explicit description of the liftings given by Theorem \ref{thm:2steps} in this particular case:

\begin{definition}\label{def:generic}
For each $\bs\mu\in\mR$ as in \eqref{eqn:mu}, $\u(\bs\mu)$ is the algebra quotient of $T(V)\#H$ given by the following set of relations:

When $(N,3)=1$,
\begin{align*}
[a_{12},a_2]_c &= 0; \qquad a_{11112}=0; \qquad a_i\a{N}=\mu_i(1-g_i^N), \ i=1,2; \\
\qquad a_{12}^N &= \mu_{12}(1-g_{12}^N)-\ta_1\mu_2a_1^Ng_2^N; \\
a_{112}^N &=\mu_{112}(1-g_{112}^N)- \ta_2\mu_{12}a_1^Ng_{12}^N- \ta_3\mu_2a_1^{2N}g_2^N;\\
a_{1112}^{N}&= \mu_{1112}(1-g_{1112}^N)- \ta_4\mu_{112}a_1^Ng_{112}^{N}- \ta_5\mu_{12}a_1^{2N}g_{12}^{N}-\ta_6\mu_2a_1^{3N}g_{2}^{N};\\
a_{\beta}^{N}&=\mu_{\beta}(1-g_{\beta}^N)-
\ta_7\mu_{12}a_{112}^Ng_{12}^N -  \ta_8\mu_2 a_{1112}^{N}g_2^{N}- \ta_9\mu_{12}^2a_1^{N}g_{12}^{2N}\\
&\quad - \ta_{10}\mu_2\mu_{12}a_1^{2N}g_{122}^{N}
- \ta_{11}\mu_2 a_{112}^Na_1^Ng_2^{N}- \ta_{12}\mu_2^2a_1^{3N}g_{2}^{2N}.
\end{align*}

When $(N,3)=3$, $N=3M$,
\begin{align*}
[a_{12},a_2]_c &= 0; \qquad a_{11112}=0; \qquad a_1\a{N}=\mu_1(1-g_1^N); \qquad a_2\a{M}=\mu_2(1-g_2^M);\\
a_{12}^N &=\mu_{12}(1-g_{12}^N)-  \tb_1\mu_2 a_{\beta}^{M}g_2^{M}- \tb_2\mu_2^2a_{1112}^{M}g_2^{2M}- \tb_3\mu_2^3a_1^Ng_2^N\\
a_{112}^N &=\mu_{112}(1-g_{112}^N)-  \tb_4\mu_{12}a_1^Ng_{12}^N- \tb_5\mu_{\beta}a_{1112}^{M}g_{\beta}^{M} - \tb_6\mu_2^3a_1^{2N}g_2^N\\
 &\quad - \tb_7\mu_2 a_{1112}^{2M}g_{2}^{M}- \tb_8\mu_2^2a_{1112}^{M}a_1^Ng_{2}^{2M}
- \tb_9\mu_2\mu_{\beta}a_1^Ng_{12}^N;\\
a_{1112}^{M}&= \mu_{1112}(1-g_{1112}^M)- \tb_{10}\mu_2a_1^Ng_2^{M};\\
a_{\beta}^{M}&=\mu_{\beta}(1-g_{\beta}^M)- \tb_{11}\mu_2^2a_1^Ng_2^{2M}-\tb_{12}\mu_2a_{1112}^{M}g_2^{M}.
\end{align*}
\end{definition}


\subsection{The degenerate case}\label{subsec:degeneratecase}\label{sec:deformed}
We now investigate the case of the braiding matrix $\qb_d$ as in \eqref{eqn:quantum-def}.

Recall the strategy from \S \ref{sec:strategy}. For each $\bs\lambda=(\lambda_{1},\lambda_{2})\in\k^2$ subject to 
\begin{align}\label{eqn:lambda}
\lambda_{1}&=0, & \text{ if } & g_{11112}=1  \text{ or }\chi_{11112}\neq \eps,\\
\notag \lambda_{2}&=0, & \text{ if } & g_{122}=1   \text{ or }\chi_{122}\neq \eps,
\end{align}
see \eqref{eqn:R}, we consider the algebra \cf Example \ref{exa:step0}:
\begin{align*}
\mE(\bs\lambda)&=\k\lg y_1,y_2 | y_{11112}-\lambda_{1},[y_{12},y_2]_c-\lambda_{2}\rg.
\end{align*}
\begin{lem}\label{lem:neqzero}
For every $\bs\lambda$ as in \eqref{eqn:lambda}, $\mE(\bs\lambda)\neq 0$.
\end{lem}
\pf
We check this with \texttt{GAP}, see log files at \newline \texttt{http://www.mate.uncor.edu/$\sim$aigarcia/publicaciones.htm}.
\epf
We set $\mA(\bs\lambda):=\mE(\bs\lambda)\# H$ and 
\begin{align*}
\u(\bs\lambda)&=\k\lg a_1,a_2|a_{11112}-\lambda_{1}(1-g_1^4g_2),[a_{12},a_2]_c-\lambda_{2}(1-g_1g_2^2)\rg.
\end{align*}
Then $\u(\bs\lambda)$ is a cocycle deformation of $\mH'$ and $\mA(\bs\lambda)$ is a $(\u(\bs\lambda),\mH')$ bi-cleft object, by Proposition \ref{pro:strategy}.
We fix $\bs\lambda$ as in \eqref{eqn:lambda} and set $\mA=\mA(\bs\lambda)$, $\u=\u(\bs\lambda)$. 
We denote the corresponding section and coactions, respectively, by
\begin{align*}
\gamma&:\mH'\to\mA, &\rho&:\mA\to \mA\ot \mH', & \delta&:\mA\to \u\ot \mA
\end{align*}
For each   $\bs\mu\in\mR$ as in \eqref{eqn:mu}, we set
\begin{align*}
\mE(\bs\lambda,\bs\mu)&=\mE(\bs\lambda)/\lg \gamma(r_\alpha)-\mu_\alpha:\alpha\in\Delta_+\rg, \qquad \mA(\bs\lambda,\bs\mu):=\mE(\bs\lambda,\bs\mu)\# H;
\end{align*}
let $\tau\colon\mA(\bs\lambda)\twoheadrightarrow\mA(\bs\lambda,\bs\mu)$ denote the canonical projection. Let $\tilde r_\alpha\in\u(\bs\lambda)$, $\alpha\in\Delta_+$, see \eqref{eqn:solution2}, be the unique solution of the equation
\begin{equation}\label{eqn:solution2-g2}
\tilde r_\alpha\ot 1=(\id\ot\tau)\delta(\gamma(r_\alpha))-g_\alpha^{N_{\alpha}}\ot \gamma(r_\alpha).
\end{equation}
and consider the quotient algebras
\[
\u(\bs\lambda,\bs\mu)=\u(\bs\lambda)/\lg \tilde r_{\alpha}-\mu_{\alpha}(1-g_\alpha^{N_{\alpha}}):\alpha\in\Delta_+\rg.
\] 
\begin{theorem}
Fix $\bs\lambda$ as in \eqref{eqn:lambda} and  $\bs\mu\in\mR$ as in \eqref{eqn:mu}. Then
\begin{enumerate}
\item $\u(\bs\lambda,\bs\mu)$ is a cocycle deformation of $\mH$ and $\gr \u(\bs\lambda,\bs\mu)\simeq \mH$.
\item $\mA(\bs\lambda,\bs\mu)$ is a $(\u(\bs\lambda,\bs\mu),\mH)$ bi-cleft object.
\end{enumerate}
Conversely, if $L$ is a lifting of $V\in\ydh$, then there are $\bs\lambda$ as in \eqref{eqn:lambda} and  $\bs\mu\in\mR$ as in \eqref{eqn:mu} such that $L\simeq \u(\bs\lambda,\bs\mu)$.
\end{theorem}
\pf
This is Theorem \ref{thm:strategy}. 
\epf

As explained in \S \ref{sec:strategy}, see Example \ref{exa:step0}, to give an explicit description of the cleft objects $\mA(\bs\lambda,\bs\mu)$ and the liftings $\u(\bs\lambda,\bs\mu)$, we need to compute $\gamma(r_\alpha)$, $\alpha\in\Delta_+$ and solve the equations \eqref{eqn:solution2-g2}.

This is the content of Proposition \ref{pro:explicit-cleft} next and Theorem \ref{thm:deformed} below.

\begin{pro}\label{pro:explicit-cleft}
Let $\bs\lambda$ be as in \eqref{eqn:lambda} and $\bs\mu$ be as in \eqref{eqn:mu}. Then
$\mA(\bs\lambda,\bs\mu)$ is the quotient of $T(V)\# H$ by the ideal generated by the relations
\begin{align*}
y_{11112}&=\lambda_{1}; \ \ [y_{12},y_2]_c=\lambda_{2}; \ \ y_i^7=\mu_i, i=1,2; \ \ y_{12}^7=\mu_{12}; \ \ y_{1112}^7=\mu_{112}; \\
y_{112}^7&=\mu_{112} - (8q+9q^2+3q^3+4q^4+5q^5-q^6)\lambda_{1}^2y_{12}^3y_\beta\\
&\quad  - (5q+8q^2+9q^3+8q^4+5q^5)\lambda_{1}^2y_{2}y_\beta^2 + 7q^3(1+2q+q^2)\lambda_{1}^3y_{2}^2y_{12}^2 \\
&\quad - (2q+4q^2+13q^3+22q^4+17q^5+5q^6)\lambda_{2}\lambda_{1}^3y_2y_{12},\\
y_\beta^7&=\mu_\beta + (210q+224q^2+238q^3+105q^4+462q^5+231q^6)\lambda_{2}^2\lambda_{1}^3y_{12}^7 \\
&\quad + (833q+1911q^2+2352q^3+1911q^4+833q^5)\lambda_{2}\lambda_{1}^5y_{2}^7.
\end{align*} 
\end{pro}
\pf
For each $\alpha\in \Delta_+$, $\gamma(r_\alpha)\in\mA(\bs\lambda)$ is the unique element satisfying
\begin{equation*}
\rho(\gamma(r_\alpha))=\left(\gamma\ot\id\right)\Delta_{\mH'}(r_\alpha),
\end{equation*}
see \eqref{eqn:solution2}. Recall the notation $\ta_i\in\k$, $i\in\I_{12}$, from p.~\pageref{page:ai}. It readily follows that $\gamma(x_i^7)= y_i^7$ and 
$\gamma(x_{12}^7)=y_{12}^7$, as 
\begin{align*}
\rho(y_i^7) &= y_i^7\ot 1+g_i^7\ot x_i^7, \quad i=1,2;\\
\rho(y_{12}^7) &= y_{12}^7\ot 1 + g_{12}^7\ot x_{12}^7 + \ta_1y_{1}^7g_{2}^7\ot x_{2}^7.
\end{align*}
Next, we use \texttt{GAP} to check that
\begin{align*}
\rho(y_{1112}^7) &= y_{1112}^7\ot 1+ g_{1112}^7\ot x_{1112}^7 + \ta_4y_{1}^7g_{112}^7\ot x_{112}^7 + \ta_5y_{1}^{14}g_{12}^7\ot x_{12}^7 \\
&\quad + \ta_6y_{1}^{21}g_{2}^7\ot x_{2}^7 .
\end{align*}
Hence $\gamma(x_{1112}^7)=y_{1112}^7$.
On the other hand, we checked using \texttt{GAP} that
\begin{align*}
\gamma(x_{112}^7) &= y_{112}^7 + (8q+9q^2+3q^3+4q^4+5q^5-q^6)\lambda_{1}^2y_{12}^3y_\beta \\
&+ (5q+8q^2+9q^3+8q^4+5q^5)\lambda_{1}^2y_{2}y_\beta^2 - 7q^3(1+2q+q^2)\lambda_{1}^3y_{2}^2y_{12}^2 \\
&+ (2q+4q^2+13q^3+22q^4+17q^5+5q^6)\lambda_{2}\lambda_{1}^3y_2y_{12}.
\end{align*} 
Similarly, we use  \texttt{GAP} to compute:
\begin{align*}
\rho(y_{\beta}^7) &=  y_{\beta}^7\ot 1+g_{\beta}^{7}\ot x_{\beta}^7 + \ta_7\gamma(x_{112}^7)g_{12}^7\ot x_{12}^7 + \ta_8y_{1112}^7g_{2}^7\ot x_{2}^7 \\
&\quad +\ta_9y_{1}^7g_{12}^{14}\ot x_{12}^{14} + \ta_{10}y_{1}^{14}g_{122}^{7}\ot x_{2}^7x_{12}^7 + \ta_{11}\gamma(x_{112}^7)y_{1}^7g_{2}^7\ot x_{2}^7\\
&\quad + \ta_{12}y_{1}^{21}g_{2}^{14}\ot x_{2}^{14} \\
&\quad + (210q+224q^2+238q^3+105q^4+462q^5+231q^6)\lambda_{2}^2\lambda_{1}^3g_{12}^7\ot x_{12}^7 \\
&\quad + (210q+224q^2+238q^3+105q^4+462q^5+231q^6)\ta_1\lambda_{2}^2\lambda_{1}^3y_{1}^7g_{2}^7\ot x_{2}^7 \\
&\quad + (833q+1911q^2+2352q^3+1911q^4+833q^5)\lambda_{2}\lambda_{1}^5g_{2}^7\ot x_{2}^7.
\end{align*} 
and thus we see that
\begin{align*}
\gamma(x_{\beta}^7) &= y_{\beta}^7 - (210q+224q^2+238q^3+105q^4+462q^5+231q^6)\lambda_{2}^2\lambda_{1}^3y_{12}^7 \\
&\quad - (833q+1911q^2+2352q^3+1911q^4+833q^5)\lambda_{2}\lambda_{1}^5y_{2}^7.
\end{align*} 
Hence the proposition follows.
\epf

\begin{definition}\label{def:deformed}
For each $\bs\mu$ as in \eqref{eqn:mu} and $\bs\lambda$ as in \eqref{eqn:lambda}, we set $\u(\bs\lambda,\bs\mu)$ as the algebra  quotient of $T(V)\# H$ by the ideal generated by the relations:
{\footnotesize

}
\end{definition}

\begin{theorem}\label{thm:deformed}
Let $V$ be a braided vector of diagonal type with braiding matrix $\qb_d$ in \eqref{eqn:quantum-def}.
The algebras $\u(\bs\lambda,\bs\mu)$ are Hopf algebras and cocycle deformations of $\B(V)\# H$. For each 
$\bs\lambda$ as in \eqref{eqn:lambda} and $\bs\mu$ as in \eqref{eqn:mu}, $\u(\bs\lambda,\bs\mu)$ 
is a lifting of $V$. 
Conversely, if $L$ is a lifting of $V$, then there are $\bs\mu$ as in \eqref{eqn:mu} and $\bs\lambda$ as in \eqref{eqn:lambda}  such that $L\simeq \u(\bs\lambda,\bs\mu)$.
\end{theorem}
\pf
We use \texttt{GAP} to compute $\delta(\gamma(x_\alpha^{N_{\alpha}}))$ for each $\alpha\in\Delta_+$, thus arriving to the algebras $\u(\bs\lambda,\bs\mu)$ from Definition \ref{def:deformed}. See the Appendix for more precisions about the use of \texttt{GAP}.
\epf

\section*{Appendix}

In this section we comment on some of the ideas behind the use of the computer for doing the calculations of \S \ref{subsec:degeneratecase}.

We stress the fact that the use of a computational resource is at the moment  unavoidable to find the deformations of the powers of the root vectors when the quantum Serre relations have been deformed. That is, it becomes too intricate to find an expression of $\rho(x_\alpha^n)$ or $ (\delta\circ \gamma)(x_\alpha^n)$, with $n\leq 7$. We choose the computer program \texttt{GAP}, together with the package \texttt{GBNP}, to perform these computations.
%

A first thing to take into account when working with \texttt{GAP} is that for the files to compile we sometimes need to describe an algebra using more relations than those given by a minimal presentation. That is, besides the (deformed or not) quantum Serre relations and those that define the root vectors $x_\alpha$, we include the relations given by the q-commutators $[x_\alpha, x_{\alpha'}]_c$. For this motive, we need the relations in Lemma \ref{lem:convex} and the following lemmata. Once we have all possible relations, we use at each step a part of them, according to the need.

%
%


\begin{lem}\label{lem:relationscleft}
In $\mE(\lambda_{1},\lambda_{2})$ it holds:
\begin{align*}
[y_{12},y_2]_c &= \lambda_{2};\quad [y_{\beta},y_{12}]_c = (q^3-1)^2\lambda_{1}y_2^2; \\
[y_{112}, y_{\beta}]_c &= (-3q-2q^2-2q^3-q^4+q^6)\lambda_{1}\lambda_{2} + q^5(q+1)(q-1)^2\lambda_{1}y_2y_{12}; \\
[y_{1112},y_{112}]_c &= (1-q^4)\lambda_{1}y_{12} ;\quad [y_{1},y_{1112}]_c = \lambda_{1};\quad [y_1,y_2]_c = y_{12}; \\
[y_{1},y_{12}]_c &= y_{112}; \quad [y_{1},y_{112}]_c = y_{1112};\quad [y_{112},y_{12}]_c = y_{\beta}; \\
[y_\beta,y_2]_c &= q^2(q+1)(q-1)^2y_{12}^3;\quad [y_{112},y_2]_c = q^4(q^2-1)y_{12}^2; \\
[y_{1112},y_{12}]_c &= (q^2+1)(1-q)y_{112}^2 + q^4(1-q)(1+q^2+q^4)\lambda_{1}y_2; \\
[y_{1},y_{\beta}]_c &= (q^2+1)(1-q)y_{112}^2 + q^4(1-q)(1+q^2+q^4)\lambda_{1}y_2; \\
[y_{1112},y_{\beta}]_c &= (q-1)(1-q^3)(1+q^3+q^5)y_{112}^3 + (q^4-q^3)^2\lambda_{1}y_{12}^2\\
&\quad - q^4(q-1)^2(q^2+1)\lambda_{1}y_2y_{112}; \\
[y_{1112},y_{2}]_c &= q^4(q^2-q-1)y_{\beta} +q(q^3-1)y_{12}y_{112}.
\end{align*}
\end{lem}
\pf
This is a straightforward computation, using  the q-Jacobi identity, more precisely the image of this equation under $T(V)\twoheadrightarrow \mE$.
\epf

\begin{lem}
In $\u(\bs\lambda)$ it holds:
\begin{align*}
[a_1,a_2]_c &= a_{12}; \ \ [a_{1},a_{12}]_c = a_{112}; \ \ [a_{1},a_{112}]_c = a_{1112}; \ \ [a_{112},a_{12}]_c = a_{\beta}; \\
[a_{12},a_2]_c &= \lambda_{2}(1-g_{122}); \quad [a_{1},a_{1112}]_c = \lambda_{1}(1-g_{11112}); \\
[a_{\beta},a_{12}]_c &= \lambda_{2}(q^4-1)a_{1112}g_{122} + \lambda_{1}(q^3-1)^2a_2^2 + \lambda_{2}(q^4-q)^2a_{112}a_1g_{122}; \\
[a_{112}, a_{\beta}]_c &= \lambda_{2}(-2q-3q^2-4q^3-3q^4-2q^5)a_{1112}a_1g_{122} \\
&\quad +\lambda_{1}(2q+q^2+q^3+q^4+2q^5)a_2a_{12} \\
&\quad +\lambda_{2}\lambda_{1}(2q+2q^2+q^3+3q^4+4q^5+2q^6)g_{122}g_{11112} \\
&\quad +\lambda_{2}\lambda_{1}(q+q^3-2q^4-4q^5-3q^6)g_{122} \\
&\quad +\lambda_{2}\lambda_{1}(-3q-2q^2-2q^3-q^4+q^6); \\
[a_{1112},a_{112}]_c &= \lambda_{1}(1-q^4)a_{12}
+\lambda_{2}(-1+3q^2+q^3+q^4+3q^5)a_1^4g_{122}; \\
[a_\beta,a_2]_c &= q^2(q^2-1)(q-1)a_{12}^3+\lambda_{2}(q-q^4)a_{112}g_{122} \\
&\quad +\lambda_{2}q^4(q^4-1)^2a_{12}a_{1}g_{122}; \\
[a_{112},a_2]_c &= q^4(q^2-1)a_{12}^2
+\lambda_{2}(q^3-1)a_{1}g_{122}; \\
[a_{1112},a_{12}]_c &= (q^2+1)(1-q)a_{112}^2 + \lambda_{1}q^4(1-q)(1+q^2+q^4)a_2 \\
&\quad +\lambda_{2}(2q-2q^2-q^4+q^6)a_{1}^3g_{122}; \\
[a_{1},a_{\beta}]_c &= (q^2+1)(1-q)a_{112}^2 + \lambda_{1}q^4(1-q)(1+q^2+q^4)a_2 \\
&\quad +\lambda_{2}(2q-2q^2-q^4+q^6)a_{1}^3g_{122}; \\
[a_{1112},a_{\beta}]_c &= (q-1)(1-q^3)(1+q^3+q^5)a_{112}^3  -\lambda_{1}q^4(q-1)^2(q^2+1)a_2a_{112} \\
&\quad +\lambda_{1}(q^4-q^3)^2a_{12}^2  +\lambda_{2}(-2q^2+5q^4+5q^5+6q^6)a_{1112}a_{1}^2g_{122} \\
&\quad +\lambda_{2}(-3q-5q^3-4q^4-4q^5-5q^6)a_{112}a_{1}^3g_{122} \\
&\quad +\lambda_{2}\lambda_{1}(1+4q+3q^2-q^3)a_{1}g_{122}g_{11112} \\
&\quad +\lambda_{2}\lambda_{1}( 4q-3q^2-4q^3-5q^4-7q^5+q^6)a_{1}g_{122}; \\
[a_{1112},a_{2}]_c &= q^4(q^2-q-1)a_{\beta} +q(q^3-1)a_{12}a_{112} +\lambda_{2}q^4(q^3-1)a_{1}^2g_{122}.
\end{align*}
\end{lem}
\pf
Follows as Lemma \ref{lem:relationscleft}.
\epf

Now, in order to calculate $ \rho $, we work with the algebra $ \mA\ot \mH'$. That is, we take the algebra generated by letters: $ y_2, y_{12}, y_{\beta}, y_{112}, y_{1112}, y_1, g_{1}, g_{2}, x_2, x_{12}$, $x_{\beta}, x_{112}, x_{1112} $ and $ x_1 $, subject to the relations from Lemma \ref{lem:relationscleft} and the equations (\ref{eqn:bracket-2})$-$(\ref{eqn:bracket-3}). We also add the $ q_{\alpha\alpha'}-$commutations given by the braiding and, finally, the commutations of the first $ 8 $ generators with the last $ 6 $ ones. Analogously, we consider the algebra $ \u\ot\mA$ to calculate $ \delta\circ\gamma $. This has also been the strategy followed in \cite{GV}.
See the \texttt{GAP} files at:   \texttt{http://www.mate.uncor.edu/$\sim$aigarcia/publicaciones.htm}.

We remark that these \texttt{GAP} files have not been used as presented in this reference, as they may not compile, but rather have been split into pieces that do compile, according to the need, as explained above.

In our setting, however, we encountered with two new phenomena. On the one hand, \texttt{GAP} would not deal properly with some long terms such as  $ y_{12}^3 $ and $ y_{112}^3 $ that appear in the relations of $ \mA $ and several more others of $ \u $. To solve this, we introduced mute variables of lesser length to replace them. These mute variables are clearly presented in the \texttt{GAP} files. Hence, to complete the computation, we performed a big number of substitutions by hand at the end of each process. 

On a different hand, the sizes of the files involved during the process were extremely big. The intermediary results were very large before they decreased to form the final result. Hence at some points of the process we divided the data (by hand) and performed parallel computations. In the end, the results were summed up together.

\end{document}